\documentclass[final,10pt]{elsarticle}
\usepackage{amsbsy,amsfonts,amssymb,multirow}
\usepackage{color}
\usepackage{amsmath,fontenc,epsfig,subfigure,lettrine,stmaryrd,enumerate}
\usepackage[latin1]{inputenc}
\usepackage{graphics}
\usepackage{psfrag}
\usepackage[crop=pdfcrop]{pstool}
\usepackage[margin=8pt,font=small,labelfont=bf,format=plain,indention=.5cm,labelsep=quad,skip=6pt]{caption}
\usepackage{url}
\usepackage{float}
\usepackage{algorithm}
\usepackage{algorithmic}
\usepackage{fullpage}

 \usepackage{amsmath,amssymb,epsfig}
  \usepackage{subfigure}
 \usepackage{subfigmat}

  \def\beq{\begin{equation}}
\def\eeq{\end{equation}}
\def\esplit{\end{split}}
\def\beqalign{\begin{array}{rl}}
\def\eeqalign{\end{array}}

\def\Abold{\mathbf{A}}
\def\Bbold{\mathbf{B}}

\def\Dbold{\mathbf{D}}

\def\Gbold{\mathbf{G}}

\def\Ibold{\mathbf{I}}

\def\Pbold{\mathbf{P}}
\def\Qbold{\mathbf{Q}}

\def\bbold{\mathbf{b}}

\def\ebold{\mathbf{e}}

\def\hbold{\mathbf{h}}

\def\mbold{\mathbf{m}}
\def\nbold{\mathbf{n}}

\def\ubold{\mathbf{u}}
\def\vbold{\mathbf{v}}
\def\wbold{\mathbf{w}}

\def\0bold{\boldsymbol{0}}
\def\0{\mathbf{\0}}

\begin{document}
\begin{frontmatter}


\title{High-order accurate difference schemes for the Hodgkin-Huxley equations}

\author[su]{David Amsallem\corref{cor1}}
\ead{amsallem@stanford.edu}
\ead[url]{stanford.edu/~amsallem}

\author[liu]{Jan Nordstr\"om}
\ead{jan.nordstrom@liu.se}


\address[su]{Department of Aeronautics and Astronautics, Durand Building, Room 028, 496 Lomita Mall. Stanford University, Stanford, 94305-4035, USA}

\address[liu]{Department of Mathematics, Link\" oping University, SE-581 83 Link\"oping, Sweden }


\begin{abstract}
A novel approach for simulating potential propagation in neuronal branches with high accuracy is developed. The method relies on high-order accurate difference schemes using the Summation-By-Parts operators with weak boundary and interface conditions applied to the Hodgkin-Huxley equations. This work is the first demonstrating high accuracy for that equation. Several boundary conditions are considered including the non-standard one accounting for the soma presence, which is characterized by its own partial differential equation. Well-posedness for the continuous problem as well as stability of the discrete approximation is proved for all the boundary conditions. Gains in terms of CPU times are observed when high-order operators are used, demonstrating the advantage of the high-order schemes for simulating potential propagation in large neuronal trees.
\end{abstract}

\end{frontmatter}


\section{Introduction}

Understanding the integration of synaptic input by a neuron and the propagation of the signal to its own output synapses is of high importance in neurosciences. Numerical simulation of such a phenomena has become an option since Hodgkin and Huxley developed their model in 1952~\cite{hodgkin52}.  The Hodgkin-Huxley equations are a set of coupled partial and ordinary differential equations. The first one is the cable equation that describes the distribution and evolution of the intracellular potential. The other equations are related to the evolution of  gating variables describing ion channels dynamics inside a neuron, which is typically constituted of dendritic branches, a soma, an axon and synapses. Appropriate boundary conditions can be associated with branch ends, junctions and the soma. The soma boundary condition is non-standard since it consists of a linear relation between time and space derivatives of the solution and the solution itself.

Subsequently, Rall developed methods for solving the potential propagation in passive neuronal trees where branches satisfying the ``3/2" law are connected to the soma~\cite{rall60}. Hines then simulated the potential propagation in a neuronal tree without soma using second order finite difference schemes~\cite{hines84}. That scheme was subsequently used for predictions of potential propagation in dendritic trees~\cite{kellems10,amsallem11:iccad}.

In the present work, high-order accurate difference schemes based on Summation-By-Parts (SBP) operators with the weak Simultaneous Approximation
Term (SAT) procedure~\cite{carpenter99,nordstrom99,nordstrom01,mattsson04:ad,mattsson04,svard07,svard08,nordstrom09,lindstrom10,berg11,gong11} are applied to the Hodgkin-Huxley equations and their associated boundary conditions. Rigorous convergence properties are demonstrated, even in the presence of the non-standard soma boundary condition. The ability to apply high order schemes for the solution of the Hodgkin Huxley equations is essential since it results in a lower computational cost for the large systems that arise when neurons with large dendritic trees are considered. It will be demonstrated that the SBP-SAT technique provides a modular approach that is particularly effective for coupling of branches in a dendritic tree. Indeed, the SBP-SAT method offers a very effective way to enforce potential continuity and current conservation at the junction between those branches.

This paper is organized as follows: in Section~\ref{sec:ContProb}, the continuous set of partial differential equations and associated boundary conditions are presented and their strong well-posedness is demonstrated using energy estimates for two essential cases: (1) the case of an axon connected to the soma and (2) the case of a dendritic tree. In Section~\ref{sec:DiscProb}, discretization of each of the aforementioned problems is carried out and the associated penalty coefficients are chosen so that semi-discrete energy estimates hold. This results in identical estimates as in the continuous case and strong stability of the semi-discrete problem. In Section~\ref{sec:OrderConv}, the order of convergence associated with each of the  two aforementioned problems is investigated using the method of manufactured solutions. Applications of the new proposed approach to large neuronal systems are reported in Section~\ref{sec:Appli}. Conclusions are drawn in Section~\ref{sec:Concl}.

\section{The continuous problem}\label{sec:ContProb}
\subsection{Equations}
The Hodgkin-Huxley equations~\cite{hodgkin52} are a set of coupled partial and ordinary differential equations expressed in terms of the (1) intracellular potential $u(x,t)$ and (2) gating variables $m(x,t)$, $h(x,t)$ and $h(x,t)$ describing ion channels dynamics. The computational domain is in this section $x\in[0,L]$.

The equation for the potential  $u$ is based on the cable equation~\cite{keener09} and can be written as
\begin{gather}
\begin{split}\label{eq:cable}
u_t &=\frac{ \mu}{a(x)} (a(x)^2u_{x})_x - \frac{1}{C_m}g(m(x,t),h(x,t),n(x,t)) u \\
&+ \frac{1}{C_m}f(m(x,t),h(x,t),n(x,t),x,t),~~(x,t)\in[0,L]\times [0,~T]. 
\end{split}
\end{gather}
In~(\ref{eq:cable}), $a(x)$ is the radius of the neuron at the location $x$, $C_m$ the specific membrane capacitance and 
\begin{equation*}
\mu = \frac{1}{2C_m R_i } > 0,
\end{equation*} 
 where  $R_i$ denotes  the axial resistivity.  The conductance $g(m,h,n)$ of the cable in terms of the gating variables is
\begin{equation}
 g(m,h,n) = g_1 m^3 h + g_2 n^4 + g_3,
 \end{equation}
 where $g_i>0,~ i=1,2,3$. The expression for $f(m,h,n,x,t)$ is given by
 \begin{equation}
 f(m,h,n,x,t) = g_1 E_1m^3 h + g_2 E_2n^4 + g_3 E_3 - I(x,t), 
 \end{equation}
 where $E_i,~i=1,2,3$ are equilibrium potentials. $I(x,t)$ is an input current at location $x$ that originates either from artificial current injection or synaptic input from another neuronal cell.

The equations describing the evolution of the gating variables are 
\begin{gather}
\begin{split}\label{eq:gatingeqs}
 m_t(x,t) &= \alpha_m(u(x,t)) (1-m(x,t))-  \beta_m(u(x,t)) m(x,t),\\
 h_t(x,t) &=\alpha_h(u(x,t)) (1-h(x,t))-  \beta_h(u(x,t)) h(x,t),\\
 n_t(x,t) &= \alpha_n(u(x,t)) (1-n(x,t))-  \beta_n(u(x,t)) n(x,t).
\end{split}
\end{gather}
where $(x,t)\in[0,L]\times [0,T]$. Expressions for $\alpha_m$, $\alpha_h$, $\alpha_n$, $\beta_m$, $\beta_h$ and $\beta_n$  are provided in Appendix A.

One can prove the following property based on the equations associated with the gating variables.

\noindent
{\bf Proposition 1.} Let $x\in[0,L]$. If $0\leq m(x,0) \leq 1$ and if $m(x,t)$ is $\mathcal{C}^0$ in time, then $0\leq m(x,t) \leq 1~\forall t\in[0,T]$. 

Proof: see Appendix B. 

\noindent
The same property also applies for $h$ and $n$. As a result, the function $g$ can be bounded as follows
\begin{equation}\label{eq:gbounds}
0 < g_3 \leq g(m,n,h) = g_1 m^3 h + g_2 n^4 + g_3\leq \sum_{i=1}^3 g_i.
\end{equation}
The bound~(\ref{eq:gbounds}) will be useful when well-posedness of the initial boundary value problems is studied in Section~\ref{sec:well-posed}. 

Several types of boundary conditions can be associated with the cable equations~\cite{keener09}.
\subsection{Voltage clamp boundary conditions}
The simplest boundary condition associated with Eq.~(\ref{eq:cable}) is the prescribed potential condition at $x_b=0$ or $x_b = L$:
\begin{equation}
 u (x_b,t) = u_0(t).
\end{equation}
This corresponds to controlling the potential at the extreme end of the cable.

\subsection{Sealed end boundary conditions}
\begin{figure}
\begin{center}
\includegraphics[height=1.2in]{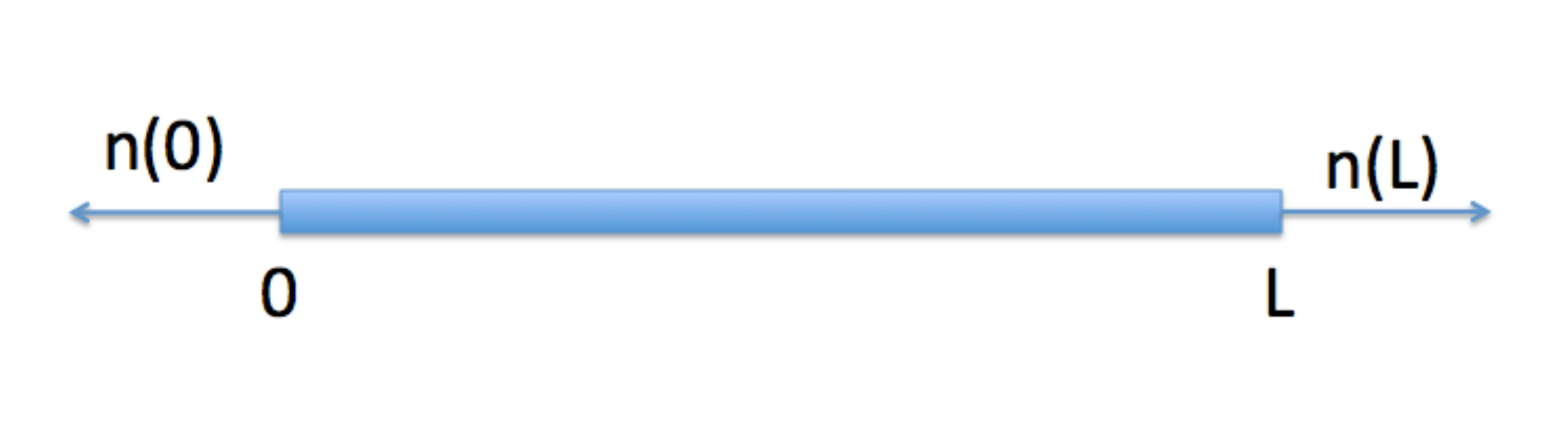}
 \caption{Outer normal vectors to a branch}
 \label{fig:OuterNormal}
 \end{center}
\end{figure}

Another possible boundary condition for Eq.~(\ref{eq:cable}) is the sealed end boundary condition at $x_b\in\{0,L\}$:
\begin{equation}
\nabla u (x_b,t)\cdot n(x_b) = 0,
\end{equation}
where $n(x_b)$ is the outer normal vector to $[0,L]$ at the end point $x_b$ of the cable, as depicted in Figure~\ref{fig:OuterNormal}. $\nabla u (x_b,t)\cdot n(x_b)$ corresponds, up to a constant factor, to the current  that exits the domain $[0,L]$ at $x_b$. Hence, enforcing the sealed end boundary conditions is equivalent to prescribing that there is no current exiting the neuron at $x_b$.

\noindent
{\bf Remark.} Although the field $u$ only depends on one space variable, in the remainder of the paper, the notation $\nabla u (x_b,t)\cdot n(x_b)$ is used in place of $-u_x(x_b,t)$ when $x_b$ is a left boundary (that is $n(x_b)=-1$) and in place of $u_x(x_b,t)$ when $x_b$ is a right boundary ($n(x_b) = 1$).

\subsection{Soma boundary conditions}
\begin{figure}
\begin{center}
 \begin{subfigmatrix}{2}
  \subfigure[Cable and soma]{\includegraphics[height=2.5in]{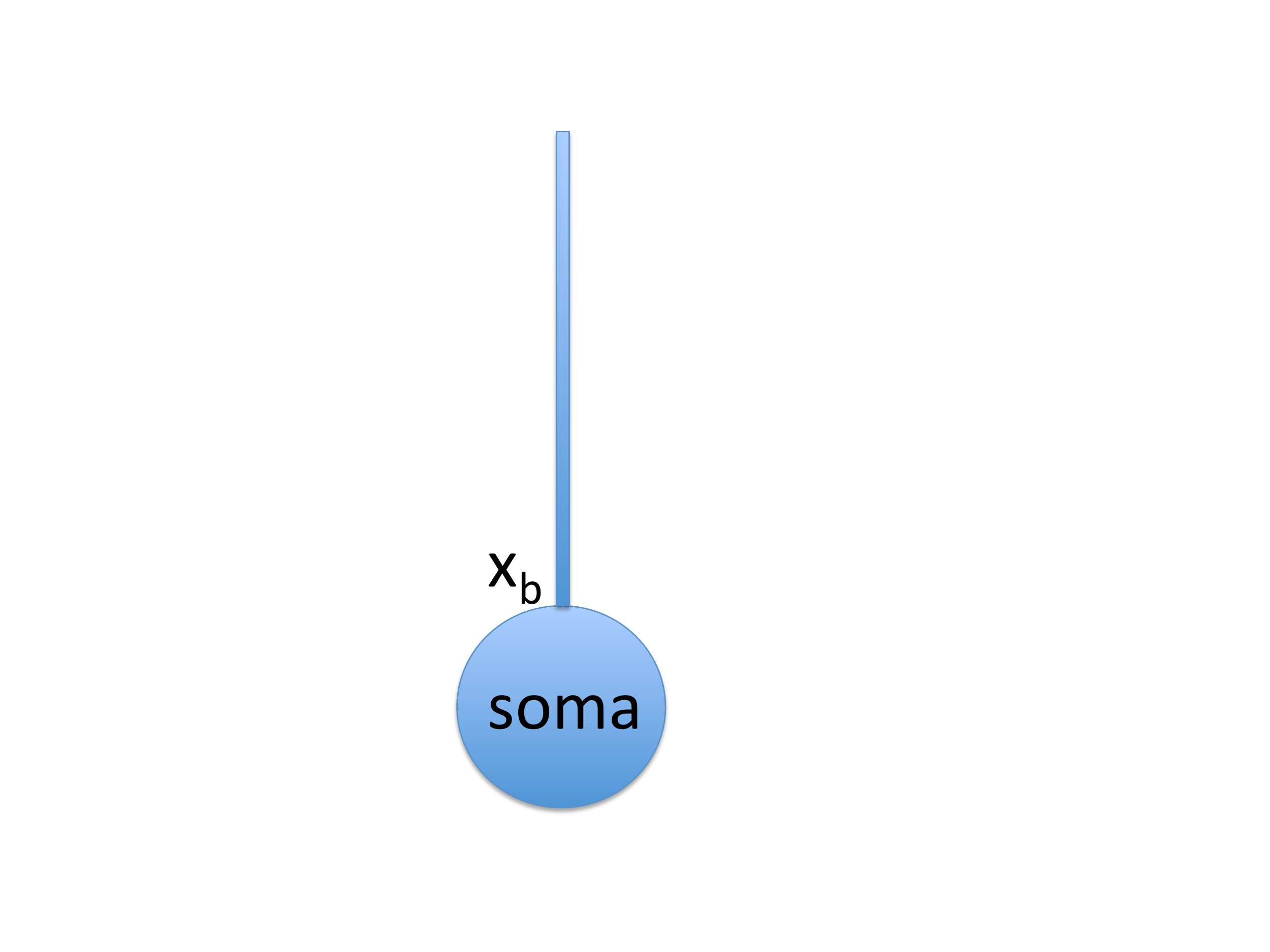}}
  \subfigure[Junction]{ \includegraphics[height=2.5in]{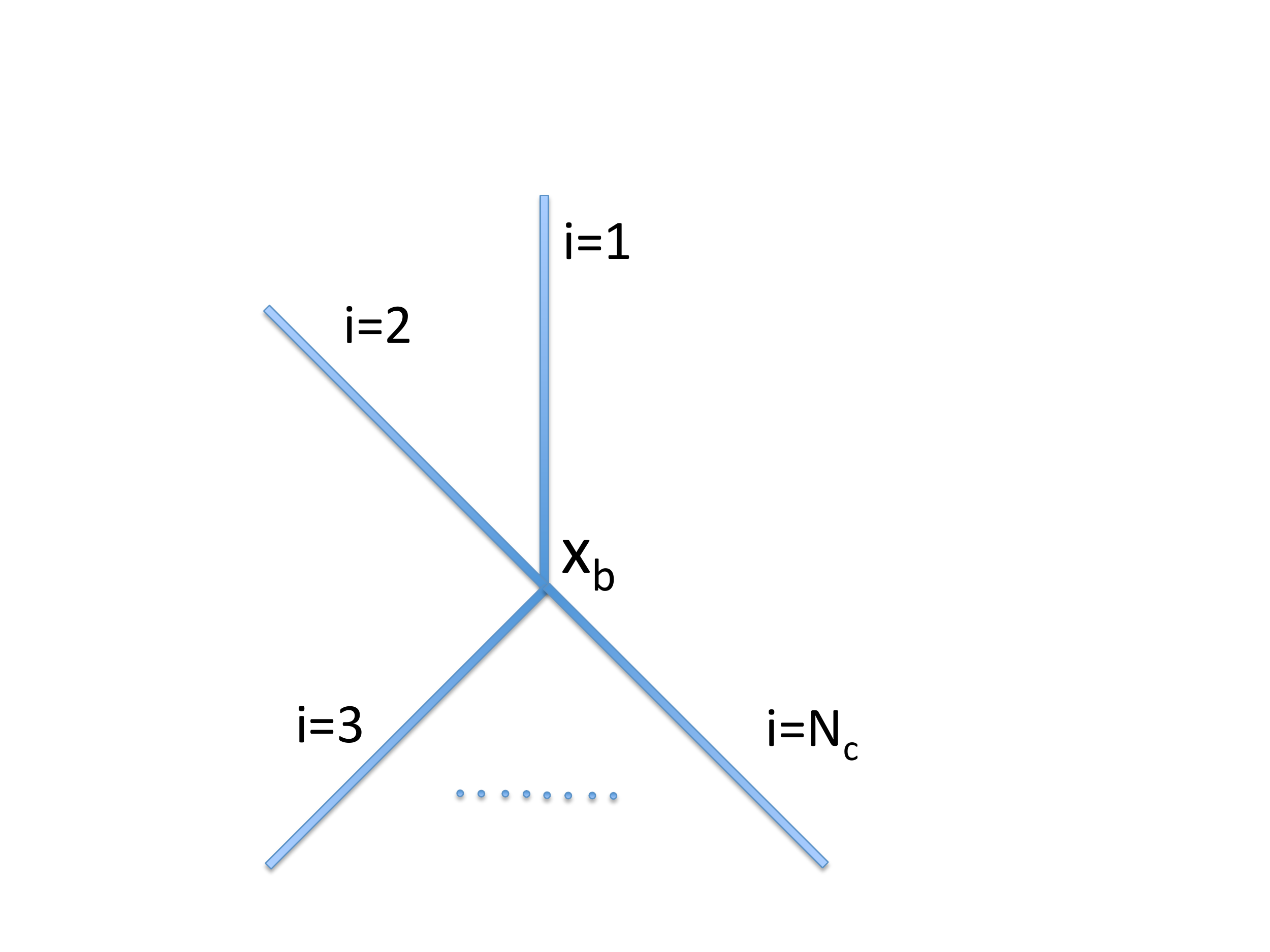}}
  \end{subfigmatrix}
 \caption{Two configurations for a neuronal network}
 \label{fig:Configurations}
 \end{center}
\end{figure}

When a soma is located at an extremity $x_b$ of a cable, a boundary condition describing the current conservation in the soma applies. This configuration is graphically depicted in Figure~\ref{fig:Configurations} (a). The boundary condition is
\begin{gather}\label{eq:somaBCeq}
\begin{split}
u_t(x_b,t) &= -\eta a(x_b)^2 \nabla u(x_b,t) \cdot n(x_b) - \frac{1}{C_m}g(m(x_b,t),n(x_b,t),h(x_b,t)) u(x_b,t) \\
&+\frac{1}{C_m} f(m(x_b,t),n(x_b,t),h(x_b,t),x_b,t)
\end{split}
\end{gather}
with 
\begin{equation}
\eta = \frac{\pi}{A_{\text{soma}}R_i C_m} > 0
\end{equation}
where $A_{\text{soma}}$ denotes the soma surface area.

\noindent
{\bf Remark.} Note here the similarity between the boundary equation~(\ref{eq:somaBCeq}) and the cable equation~(\ref{eq:cable}). The only difference lies in the spatial derivative term.


\subsection{Junction interface conditions}

When multiple cables meet at a common junction, as depicted in Figure~\ref{fig:Configurations} (b), appropriate interface conditions have to be applied. These  conditions correspond to potential continuity and current conservation~\cite{kellems10}. If $N_c$ cables of respective potential $u^{(i)},~i=1,\cdots,N_c$ and cable radius $a^{(i)}(\cdot)$ have a junction at $x=x_b$, then these interface conditions are
\begin{equation} 
u^{(i)}(x_b,t) = u^{(j)}(x_b,t),~i,j\in\{1,\cdots,N_c\}
\end{equation}
for potential continuity and
\begin{equation} 
\sum_{i=1}^{N_c} (a^{(i)}(x_b))^2 \nabla u^{(i)}(x_b,t)\cdot n^{(i)}(x_b) = 0
\end{equation}
for current conservation.



\subsection{Well-posedness}\label{sec:well-posed}
\subsubsection{Cable and soma}\label{sec:contincablesoma}
The case of a single cable connected to the soma is first studied. Dropping the source terms, the cable equation is written as
\begin{equation}\label{eq:NdCable}
 a( x) u_{ t} = \mu( a( x)^2 u_{x})_{x} - \frac{ a(  x)}{C_m}  g(  x,  t)   u,~~(  x,  t )\in[0,  L]\times [0,  T]
\end{equation} 
where $  g(  x,  t)  =   g(  m(  x,  t),  n(  x,  t),  h(  x,  t)) \in  [  g_{\min},g_{\max}]$  and $  a(  x) \in  [  a_{\min} , a_{\max}]$, with
\begin{equation}
  g_{\min} = g_3>0,~~ g_{\max} = \sum_{i=1}^3 g_i ,~~  a_{\min}  > 0. 
\end{equation}
The boundary conditions associated with the cable equation are, in the case of a single cable connected to the soma at $x_b=L$,
\begin{equation}\label{eq:cablesomaBC}
  u_{  x} (0) = 0,~~  u_{  t}(  L) = - \eta   a(  L)^2   u_{  x}(   L) -   \frac{1}{C_m}g(  L,  t)  u(  L).
\end{equation}
Defining the energy norm of $  u$ at time $  t$ as
\begin{equation}\label{eq:norm}
\|   u\|^2 = \int_0^{  L}   u(  x,  t)^2 d  x,
\end{equation}
 By multiplying Eq,~(\ref{eq:NdCable}) with $u$ and integrating over the domain $[0,  L]$, one obtains
 \begin{equation}\label{eq:IBP}
\frac{1}{2}\| \sqrt{  a}   u\|^2_{  t} =  - \|  \sqrt{\mu} a   u_{  x}\|^2 +   \mu(a(  L)^2  u_{  x}(  L)  u(  L) -   a(0)^2  u_{  x}(0)  u(0) )-\left \| \sqrt{  \frac{g   a}{C_m}}   u\right\|^2.
\end{equation}
The sealed boundary condition at $x=0$ and the soma boundary equation at $x=L$ result in  
\begin{equation}\label{eq:ineq1}
\frac{1}{2}\|\sqrt{  a}  u\|^2_{  t} + \frac{\mu}{2  \eta }  (  u^2)_{  t}(L)  = - \| \sqrt{\mu} a   u_{  x}\|^2 -\frac{ \mu g(  L,  t)}{  \eta C_m}   u(  L)^2  -\left \|\sqrt{ \frac{ g  a}{C_m}}  u\right\|^2.
\end{equation}
Defining the weighted norm
\begin{equation*}
\|   u\|^2_{*}  = \| \sqrt{  a}  u \|^2 + \frac{\mu}{\eta}   u(L)^2,
\end{equation*}
 Eq.~(\ref{eq:ineq1}) becomes, 
\begin{equation}\label{eq:CEECableSoma}
(\|   u\|^2_{*} )_{  t}  + 2 \|  \sqrt{\mu}  a  u_{  x}\|^2 + 2\left\| \sqrt{ \frac{ g}{C_m}}  u\right\|^2_{*} = 0 .
\end{equation}
This result can be summarized in the following Proposition.

\noindent
{\bf Proposition 2.} The initial boundary value problem~(\ref{eq:NdCable})-(\ref{eq:cablesomaBC})  is well-posed.

\subsubsection{Cables with junction}\label{sec:continjunction}
In this section, the case of $N_c$ cables connected at a junction is considered. For simplification, it is assumed that all the other cable extremities are either sealed or clamped. If one of those extremities was connected to the soma, a similar approach to the one presented in the previous section can be used.

$  u^{(i)}(  x^{(i)},  t),~(  x^{(i)},  t)\in [0,  L^{(i)}]\times [0,  T]$ here denotes the potential in the $i$-th cable, $i=1,\cdots, N_c$. To simplify the notations, it is assumed here that the junction is located at $  x^{(i)}=0$, all the other cable extremities being sealed, except the last one where the potential is held constant at zero. Hence, the PDEs of interest are of the same form as~(\ref{eq:NdCable}) where $u$, $x$, $a$ and $L$ are replaced by $u^{(i)}$, $x^{(i)}$, $a^{(i)}$ and $L^{(i)}$ respectively. 
The interface conditions are
\begin{equation}\label{eq:BCjunctionIBVP}
  u^{(i)}(0,  t) =   u^{(j)}(0,  t),~\forall i,j=1,\cdots,N_c,~~\sum_{i=1}^{N_c}  (  a^{(i)}(0))^2    u^{(i)}_{  x^{(i)}} (0,t) = 0,
\end{equation}
\begin{equation}\label{eq:BCSealedIBVP}
u^{(i)}_{  x^{(i)}}(  L^{(i)},  t) = 0,~~i=1,\dots,N_c-1
\end{equation}
and
\begin{equation}\label{eq:BCClampIBVP}
u^{(N_c)}(  L^{(N_c)},  t) = 0.
\end{equation}
The energy norm for $  u^{(i)}$ at time $  t$ is defined in (\ref{eq:norm}) where the integral bounds are $0$ and $  L^{(i)}$.

For each cable $i\in 1,\cdots, N_c$ equation~(\ref{eq:IBP}) holds with appropriate variables.
%
Summing those equations and applying the sealed end and voltage clamp boundary conditions leads to
\begin{gather}
\begin{split}\label{eq:NdMultEnergySum}
\frac{1}{2}\sum_{i=1}^{N_c}\left\|  \sqrt{  a^{(i)}}   u^{(i)}\right\|^2_{  t} &=  -\mu\sum_{i=1}^{N_c}\left   \|   a^{(i)}   u^{(i)}_{  x^{(i)}}\right\|^2
-\mu \sum_{i=1}^{N_c}  (   a^{(i)}(0))^2  u^{(i)}_{  x^{(i)}}(0)  u^{(i)}(0) \\
&-  \sum_{i=1}^{N_c} \left\| \sqrt{ \frac{ g^{(i)}   a^{(i)}}{C_m}}   u^{(i)}\right\|^2.
\end{split}
\end{gather}
Applying the potential continuity condition at the junction, $u^{(1)}(0) = u^{(i)}(0),~i=2,\cdots,N_c$, the identity
\begin{equation}
\sum_{i=1}^{N_c}(   a^{(i)}(0))^2  u^{(i)}_{  x^{(i)}}(0)  u^{(i)}(0) =    u^{(1)}(0)
\sum_{i=1}^{N_c}  (   a^{(i)}(0))^2  u^{(i)}_{  x^{(i)}}(0)
\end{equation}
and the current conservation at the junction leads to
cancelation of the second factor in the right hand side of equality~(\ref{eq:NdMultEnergySum}). Hence

\begin{gather}
\begin{split}\label{eq:NdMultEnergySum2}
\frac{1}{2}\sum_{i=1}^{N_c} \left\| \sqrt{  a^{(i)}}   u^{(i)}\right\|^2_{  t} &=  -\mu \sum_{i=1}^{N_c} \left  \|   a^{(i)}   u^{(i)}_{  x^{(i)}}\right\|^2
 - \sum_{i=1}^{N_c} \left\| \sqrt{ \frac{ g^{(i)}   a^{(i)}}{C_m}}   u^{(i)}\right\|^2.
\end{split}
\end{gather}
One can define the weighted norm $\|\cdot\|_\star$ for $  u = [  u^{(1)},\cdots,  u^{(N_c)}]^T$ as
  \begin{equation}
  \|  u\|_\star^2 = \sum_{i=1}^{N_c}  \left\| \sqrt{   a^{(i)}}    u^{(i)}\right\|^2.
  \end{equation}
The energy equality then becomes
  \begin{equation}\label{eq:CEEJunction}
    \left(\|  u\|_{\star}^2\right)_t + 2\| \sqrt{ \mu a}  u_{  x}\|^2_\star +2 \left\| \sqrt{ \frac{ g}{C_m}}  u\right\|_\star^2  = 0.
  \end{equation}
  This result is summarized in the following Proposition.

\noindent
{\bf Proposition 3.} The initial boundary value problem~(\ref{eq:NdCable}), (\ref{eq:BCjunctionIBVP}), (\ref{eq:BCSealedIBVP}), (\ref{eq:BCClampIBVP}) holding for $N_c$ cables connected at a junction and with other cable extremities sealed or clamped is well-posed.

\section{The semi-discrete problem}\label{sec:DiscProb}

In this section, Eq.~(\ref{eq:NdCable}) is discretized on the domain $[0,  L]$ using a uniform mesh of $N+1$ points. The discrete approximation of $  u(\cdot,t)$ is 
\begin{equation*}
 \ubold(t) = [  u_0(t),\cdots,   u_N(t)]^T,~~\displaystyle{u_i(t) \approx u\left(\frac{i-1}{N}L,t\right),~i=1,\cdots, N.}
\end{equation*}
Operators of SBP form~\cite{gustafsson95} approximate the derivative of $ u(\cdot,t)$ as
\begin{equation}
 \ubold_x(t) = [  u_{x0}(t),\cdots,   u_{xN}(t)]^T\approx \Dbold_1  \ubold(t) = \Pbold^{-1}\Qbold \ubold(t),
\end{equation}
where $\Pbold$ is a diagonal symmetric positive definite matrix and $\Qbold$ satisfies $\Qbold+\Qbold^T = \Bbold =  \text{\bf diag}(-1,0,\cdots,0,1)$. The second  space derivative in this paper is approximated as
\begin{equation}
\ubold_{xx}(t) \approx \Dbold_2 \ubold(t) = \Dbold_1\Dbold_1\ubold(t) =  \left(\Pbold^{-1}\Qbold\right)^2 \ubold(t).
\end{equation}
Compact second derivatives in second order form also exist~\cite{mattsson04,mattsson12}.

To alleviate notations, the dependency of $\ubold$ and its derivatives on time is omitted in the following. A norm based on $\Pbold$ is
\begin{equation}
\|  \ubold\| = \sqrt{ \ubold^T\Pbold \ubold}.
\end{equation}
The semi-discrete version of~(\ref{eq:NdCable}) without inclusion of boundary conditions is
\begin{equation}\label{eq:SDcable}
 \Abold  \ubold_{  t} =\mu \Pbold^{-1}\Qbold( \Abold^2 \ubold_{  x}) -  \frac{1}{C_m} \Abold \Gbold(  t)  \ubold
\end{equation}
where the diagonal matrix $ \Abold = \text{\bf diag}(  a_0,\cdots,  a_N)$ and the diagonal matrix $  \Gbold(  t) = \text{\bf diag}(  g_0(t),\cdots,  g_N(t))$. 


\subsection{Cable with soma}

 The boundary conditions~(\ref{eq:cablesomaBC})  are  integrated in the formulation as penalty terms in the following way
 \begin{gather}
 \begin{split}\label{eq:SDPenalty}
  \Abold  \ubold_{  t} &=\mu \Pbold^{-1}\Qbold( \Abold^2 \ubold_{  x}) -  \frac{1}{C_m} \Abold \Gbold(  t)  \ubold  + \sigma_0 \Pbold^{-1}(   u_{  x0}-0)\ebold_0 \\
 &+ \sigma_{  L} \Pbold^{-1}\left(  u_{  t N} + \eta  a_N^2  u_{  x N} +    \frac{1}{C_m}g_N(  t)   u_{N}-0\right)\ebold_L,
 \end{split}
 \end{gather}
 where, for clarity, the data is indicated by $0$.
 
 \noindent
 Premultiplying Eq.~(\ref{eq:SDPenalty}) by $ \ubold^T\Pbold$ leads to
  \begin{gather}
 \begin{split}\label{eq:SDPenaltyEnergy}
 \ubold^T\Pbold  \Abold  \ubold_{  t} &=\mu \ubold^T\Qbold( \Abold^2 \ubold_{  x}) - \frac{1}{C_m} \ubold^T\Pbold  \Abold \Gbold(  t)  \ubold  + \sigma_0    u_{  x0}  u_0 \\
 &+ \sigma_{  L}   u_{N} \left(  u_{  t N} + \eta  a_N^2  u_{  x N} +  \frac{1}{C_m} g_N( t)   u_{N}\right)
 \end{split}
 \end{gather}
 Since $\Pbold$, $\Abold$ and $\Gbold(  t)$ are diagonal matrices, they commute and lead to
\begin{equation}
\frac{1}{C_m}  \ubold^T\Pbold  \Abold \Gbold(  t)  \ubold =   \ubold^T\sqrt{ \frac{\Gbold(  t) \Abold}{C_m}}\Pbold \sqrt{ \frac{ \Gbold(  t)\Abold}{C_m}}  \ubold = \left\| \sqrt{ \frac{ \Gbold(  t)\Abold}{C_m}} \ubold\right\|^2,
\end{equation}

\begin{equation}
 \ubold^T\Pbold  \Abold  \ubold_{  t} =  \ubold^T\sqrt{ \Abold}\Pbold \sqrt{ \Abold}  \ubold_{  t} = \frac{1}{2}\|\sqrt{ \Abold} \ubold\|_{  t}^2.
\end{equation}
The SBP properties lead to

\begin{gather}
\begin{split} 
  \ubold^T\Qbold(  \Abold^2 \ubold_{  x})  &=   \ubold^T(\Bbold -\Qbold^T)\Pbold^{-1}\Pbold( \Abold^2 \ubold_{  x})  \\
&= -  \ubold^T(\Pbold^{-1} \Qbold)^T\Pbold( \Abold^2 \ubold_{  x})  +   \ubold^T\Bbold ( \Abold^2 \ubold_{  x})\\
&= -    \ubold_{  x}^T\Pbold(  \Abold^2  \ubold_{  x}) -   a_0^2  u_{0}  u_{  x0} +   a_N^2  u_{N}  u_{  xN}\\
&= -\|  \Abold  \ubold_{  x}\|^2 +   a_N^2  u_{N}  u_{  xN} -   a_0^2  u_{0}  u_{  x0}. 
\end{split}
\end{gather}
 By the choice $\sigma_0 =   \mu a_0^2$ and $\sigma_{  L} = -\frac{\mu}{  \eta}$, Eq.~(\ref{eq:SDPenaltyEnergy}) becomes

  \begin{equation}\label{eq:SDPenaltyEnergy3}
\frac{1}{2}\|\sqrt{ \Abold} \ubold\|_{  t}^2= -\mu\|  \Abold  \ubold_{  x}\|^2 -\left\| \sqrt{ \frac{ \Gbold(  t)\Abold}{C_m}} \ubold\right\|^2  -\frac{\mu}{  \eta}  u_{N} \left(  u_{  t N} +   \frac{1}{C_m}g_N(t)   u_{N}\right).
 \end{equation}
 Similarly as in the continuous case, one can define a weighted norm
\begin{equation}
\|   \ubold\|^2_* = \left\| \sqrt{ \Abold}   \ubold\right\|^2 + \frac{\mu}{  \eta}   u_{N}^2,
\end{equation}
which leads following energy estimate and proposition
\begin{equation}\label{eq:SDEECableSoma}
 \|   \ubold\|_{*,t}^2 +2 \mu\|   \Abold \ubold_{  x}\|^2 +2 \left\| \sqrt{  \frac{\Gbold(  t)}{C_m}}  \ubold\right\|_*^2 = 0.
\end{equation}

\noindent
{\bf Proposition 4. } The SBP-SAT scheme~(\ref{eq:SDPenalty}) for solving the semi-discrete problem associated with the cable equation with soma as a boundary condition is stable~\cite{gustafsson95,nordstrom05}. 

\noindent
{\bf Remark.} Note that the energy estimate~(\ref{eq:SDEECableSoma}) is the semi-discrete analog to the continuous energy estimate~(\ref{eq:CEECableSoma}) derived in Section~\ref{sec:contincablesoma}.


\subsection{Cables with junction}

Each cable is here discretized on the domain $[0,  L^{(i)}],~i=1,\cdots,N_c$ using a uniform mesh of $N^{(i)}+1$ points. One defines the discrete approximation of $  u^{(i)}(\cdot,t),~i=1,\cdots,N_c$ as
\begin{equation}
  \ubold^{(i)}(t) = \left[  u^{(i)}_0(t),\cdots,  u^{(i)}_{N^{(i)}}(t)\right].
\end{equation}
The space derivatives of $ \ubold^{(i)}(t)$ are approximated as
\begin{equation}
  \ubold^{(i)}_x(t) \approx \left(\Pbold^{(i)}\right)^{-1}\Qbold^{(i)}   \ubold^{(i)}(t),
\end{equation}
where $(\Pbold^{(i)},\Qbold^{(i)}),~i=1,\cdots,N_c$ have similar properties as in the previous section. Diagonal matrices $ \Abold^{(i)}$, $  \Gbold^{(i)}(  t)$ are defined as in the previous section. Similarly as in the previous section, the dependency of $\ubold$ and its derivatives on time is omitted for clarity. The semi-discretized equations are all of the same form as~(\ref{eq:SDcable}).

The  interface conditions (\ref{eq:BCjunctionIBVP})--(\ref{eq:BCClampIBVP}) are added as penalty terms, resulting in 
\begin{gather}\label{eq:SDJunction}
\begin{split} 
 \Abold^{(i)}  \ubold^{(i)}_{  t} &= \mu^{(i)}\left(\Pbold^{(i)}\right)^{-1}\Qbold^{(i)}\left(\left( \Abold^{(i)}\right)^2 \ubold^{(i)}_{  x}\right) -  \frac{1}{C_m} \Abold^{(i)} \Gbold^{(i)}(  t)  \ubold^{(i)} \\
&+ \sum_{j=1,j\neq i}^{N_c} \sigma_{0,j}^{(i)}\left(\Pbold^{(i)}\right)^{-1}\left(\Dbold_{1}^{(i)}\right)^T\left (  u^{(i)}_0 -   u^{(j)}_0\right)  \ebold_0^{(i)}\\
&+ \sigma_{0,i}^{(i)}\left( \sum_{j=1}^{N_c}     (a^{(j)}_0)^2   u^{(j)}_{  x 0} - 0\right) \left(\Pbold^{(i)}\right)^{-1} \ebold_0^{(i)}\\
&+ p_L^{(i)} ,~~i=1,\cdots,N_c,
\end{split}
\end{gather}
where $p_L^{(i)}$ denotes the penalty term for the $i$-th cable associated with the boundary $x=L^{(i)}$:
\begin{equation}
p_L^{(i)}=  \sigma_{  L}^{(i)}   \left(\Pbold^{(i)}\right)^{-1}\left(u^{(i)}_{  x N^{(i)}}-0\right)  \ebold_{  L}^{(i)},~~i=1,\cdots,N_i-1 
\end{equation}
for the sealed boundary conditions $u^{(i)}_{  x N^{(i)}} = 0$, and
\begin{equation}
p_L^{(N_c)} = \sigma_{  L}^{(N_c)}   \left(\Pbold^{(N_c)}\right)^{-1}\left(\Dbold_{1}^{(N_c)}\right)^T (u^{(N_c)}-0)\ebold_{  L}^{(N_c)}
\end{equation}
for the voltage clamp condition $u^{(N_c)} = 0$. These penalty terms can be all recast in the same form 
\begin{equation}
p_L^{(i)}=  \sigma_{  L}^{(i)}   \left(\Pbold^{(i)}\right)^{-1}u^{(i)}_{  x N^{(i)}}  \ebold_{  L}^{(i)},~~i=1,\cdots,N_i. 
\end{equation}
Multiplying by $\left(  \ubold^{(i)}\right)^T\Pbold^{(i)}$ and using analog identities to the one derived in the previous section leads to
 \begin{gather}
 \begin{split}\label{eq:SDPenaltyEnergyMultCables}
\frac{1}{2}\left\|\sqrt{ \Abold^{(i)} } \ubold^{(i)} \right\|_{  t}^2&= -\mu\left\|  \Abold^{(i)}   \ubold^{(i)} _{  x}\right\|^2 -\left\| \sqrt{\frac{   \Gbold^{(i)} (  t) \Abold^{(i)}}{C_m}} \ubold^{(i)} \right\|^2 \\
&+ \sum_{j=1,j\neq i}^{N_c} \sigma_{0,j}^{(i)}  u^{(i)}_{x0}\left (  u^{(i)}_{0} -   u^{(j)}_{0}\right) \\
&+ \sigma_{0,i}^{(i)}\sum_{j=1,j\neq i}^{N_c}     (a^{(j)}_0)^2    u^{(i)}_0  u^{(j)}_{  x 0} \\
&+ (\sigma_{0,i}^{(i)} -\mu)   (a^{(i)}_0)^2    u_{0}^{(i)}  u_{  x 0}^{(i)}\\
&+ (\sigma_{  L}^{(i)}  +  \mu (a^{(i)}_{N^{(i)}})^2)  u^{(i)}_{  x N^{(i)}}  u^{(i)}_{ N^{(i)}},~~i=1,\cdots,N_c.
 \end{split}
 \end{gather}
Summing the $N_c$  equalities~(\ref{eq:SDPenaltyEnergyMultCables})  and choosing the penalty coefficients $\sigma_{  L}^{(i)} = -   \mu\left(a^{(i)}_{N^{(i)}}\right)^2$,
 \begin{equation}
\displaystyle{ \sigma_{0,i}^{(i)}  = \frac{ \mu}{N_c},~i=1,\cdots,N_c}
 \end{equation}
 and
  \begin{equation}
   \displaystyle{\sigma_{0,i}^{(j)} = \sigma_{0,i}^{(i)}  (a^{(j)}_0)^2  = \frac{ \mu}{N_c}   (a^{(j)}_0)^2    }
 \end{equation}
leads to the energy estimate
 \begin{equation}\frac{1}{2} \sum_{i=1}^{N_c} \left\|\sqrt{ \Abold^{(i)} } \ubold^{(i)} \right\|_{  t}^2= -\mu\sum_{i=1}^{N_c}\left\|  \Abold^{(i)}   \ubold^{(i)} _{  x}\right\|^2 -\sum_{i=1}^{N_c} \left\| \sqrt{\frac{   \Gbold^{(i)} (  t) \Abold^{(i)}}{C_m}} \ubold^{(i)} \right\|^2.
 \end{equation}
Defining the weighted norm
\begin{equation}
\| \ubold \|^2_\star = \sum_{i=1}^{N_c}\left\|   \sqrt{ \Abold^{(i)} } \ubold^{(i)}  \right\|^2,
\end{equation}
the following energy equality holds as well as the proposition
\begin{equation}
 \left( \| \ubold\|^2_\star\right) +2\mu\left \| \sqrt{\Abold} \ubold \right\|^2_\star +2 \left\| \sqrt{\frac{\Gbold(t)}{C_m}}\ubold\right\|_\star^2 = 0.
\end{equation}
\noindent 
 {\bf Proposition 5.}  The SBP-SAT scheme for solving the semi-discrete problem~(\ref{eq:SDJunction})  associated with the initial boundary value problem~(\ref{eq:NdCable}), (\ref{eq:BCjunctionIBVP}), (\ref{eq:BCSealedIBVP}), (\ref{eq:BCClampIBVP})  is stable. 
 
 \noindent
{\bf Remark.} Note that this energy estimate is the exact semi-discrete analog to the energy equality~(\ref{eq:CEEJunction}) derived in Section~\ref{sec:continjunction}.

\section{Order of convergence}\label{sec:OrderConv}
The method of manufactured solutions~\cite{knupp03} is used to determine the discretization order.  For the two problems mentioned above, a closed form solution is chosen and appropriate source terms are added to the initial boundary value problem of interest so that the closed form solution satisfies it. Time integration is done using the fourth-order explicit Runge-Kutta scheme and a time step of $\Delta t= 10^{-9}$ s, ensuring that the time discretization errors are negligible. The convergence rate is calculated as
\begin{equation}
\displaystyle{q = \frac{\log_{10}\left(  \frac{\| \vbold_{h_1} - \ubold \|}{\|  \vbold_{h_2} - \ubold\|} \right)}{\log_{10}\left(\frac{h_1}{h_2}\right)}}.
\end{equation}
\subsection{Cables with junction}

\begin{figure}
\begin{center}
\includegraphics[height=3.5in]{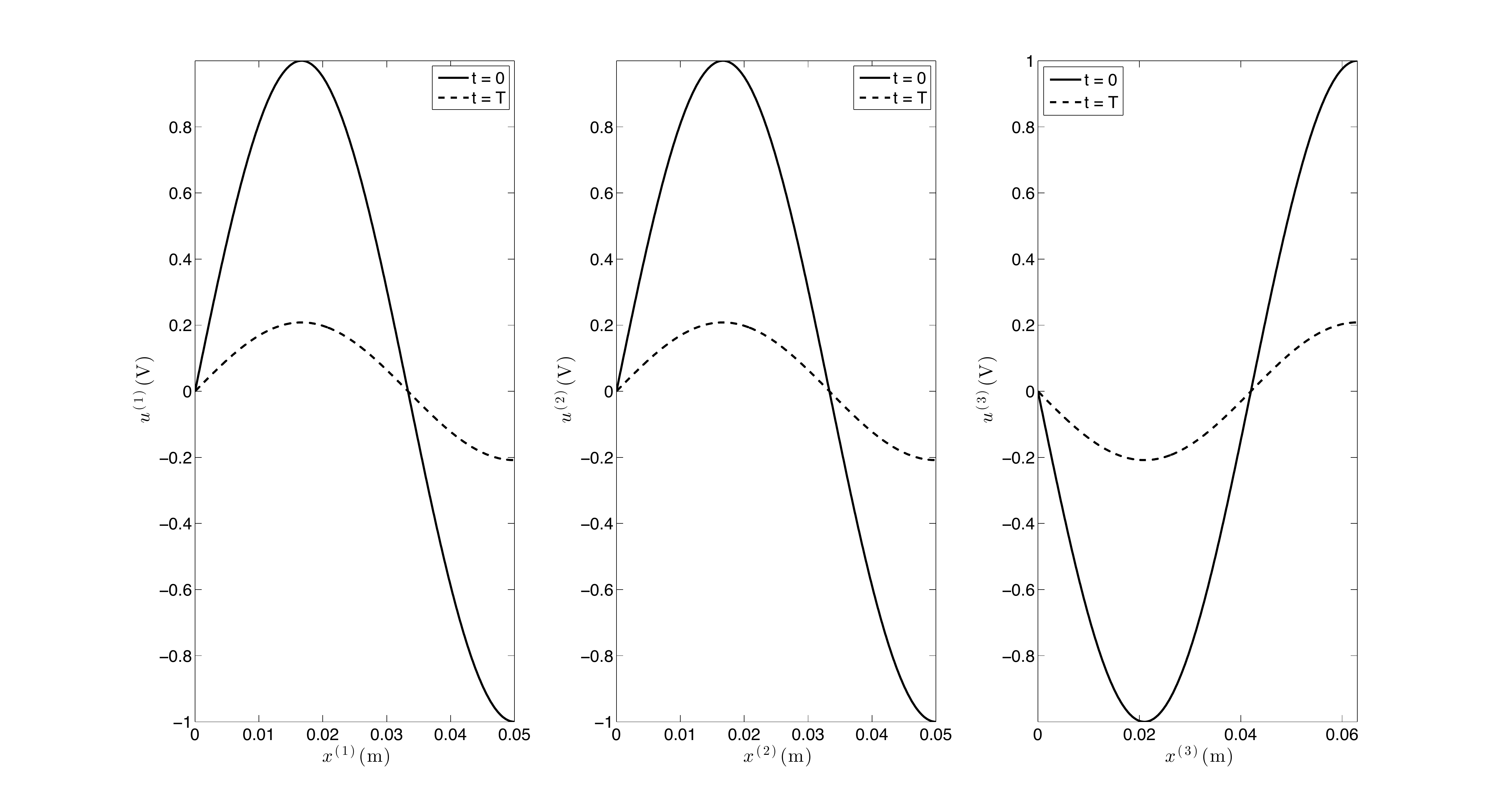}
 \caption{Exact solutions for the junction problem}
 \label{fig:exactsolJunction}
 \end{center}
\end{figure}

Three cables with junction at $x^{(i)}=0,~i=1,\cdots 3$ are considered. Their respective cable radii are constant for each cable. The exact manufactured solutions $u^{(i)}(x,t),~i=1,\cdots,3$ and their corresponding initial boundary value problem are specified in Appendix C.  The physical constants associated with the problems are reported in Tables~\ref{tab:physvar} and~\ref{tab:junctionGeom}. The exact solutions $u^{(i)}(x,t)$ at times $t=0$ and $t=T$ are plotted in Figure~\ref{fig:exactsolJunction}.

\begin{table}[htdp]
\begin{center}\begin{tabular}{|c|c|c|}
\hline
Specific membrane capacitance & $C_m$ & $10^{-2}$ F.m$^{-2}$ \\
\hline
Axial Resistivity & $R_i$ & $0.354$ ohm.m \\
\hline
Conductance & $g_1$ & $1200$ ohm$^{-1}$.m$^{-2}$ \\
\hline
Conductance &$g_2$ & $360$ ohm$^{-1}$.m$^{-2}$ \\
\hline
Conductance &$g_3$ & $3$ ohm$^{-1}$.m$^{-2}$ \\
\hline
Equilibrium potential & $E_1$ & $0.115$ V \\
\hline
Equilibrium potential &$E_2$ & $-0.012$ V \\
\hline
Equilibrium potential &$E_3$ & $0.010613$ V \\
\hline
\end{tabular} 
\end{center}\caption{Physical constants for the problems of interest}
\label{tab:physvar}
\end{table}

\begin{table}[htdp]
\begin{center}\begin{tabular}{|c|c|c|}
\hline
Cable length & $L^{(1)}$ & $0.05$ m \\
\hline
Cable length & $L^{(2)}$ & $0.05$ m \\
\hline
Cable length & $L^{(3)}$ & $0.063$ m \\ 
\hline
Cable radius & $a^{(1)}(x)$ & $0.476\times 10^{-3}$ m  \\
\hline
Cable radius & $a^{(2)}(x)$ & $0.476\times 10^{-3}$ m  \\
\hline
Cable radius & $a^{(3)}(x)$ & $0.7556\times 10^{-3}$ m  \\
\hline
\end{tabular} 
\end{center}\caption{Geometric parameters associated with the junction problem}
\label{tab:junctionGeom}
\end{table}

The solutions $u^{(i)}(x^{(i)},t) $ at $t=0$ and $t=T$ are reported in Figure~\ref{fig:exactsolJunction}.

\begin{figure}
\begin{center}
\includegraphics[height=3.0in]{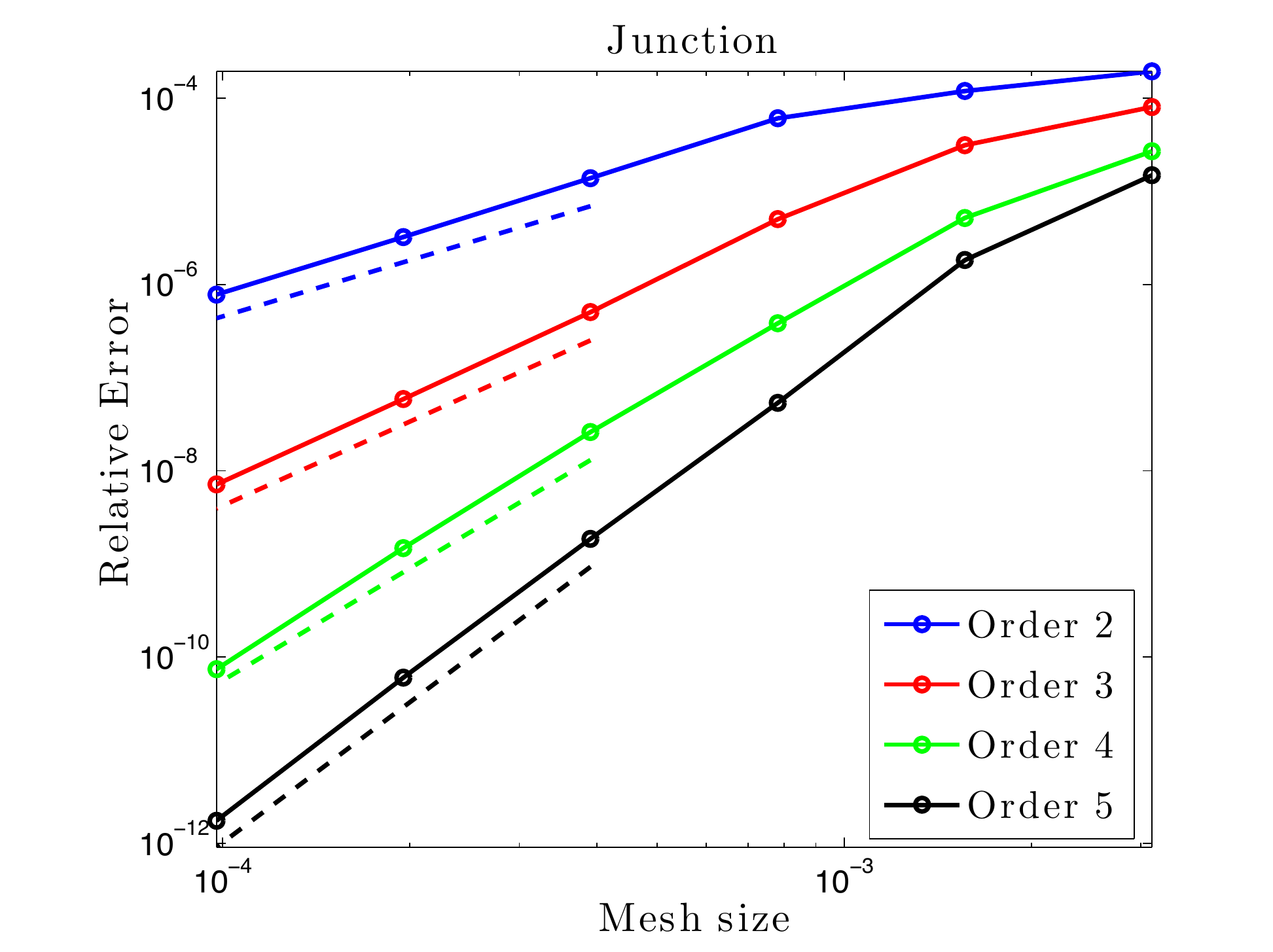}
 \caption{Errors in function of mesh size and operators order for the junction problem}
 \label{fig:JunctionCVOrderNonCompact}
 \end{center}
\end{figure}

The convergence results are reported in Table~\ref{tab:CVorderNonCompactJunction} and Figure~\ref{fig:JunctionCVOrderNonCompact}. The slopes corresponding to each operators order are also shown in Figure~\ref{fig:JunctionCVOrderNonCompact}. Second-order convergence is observed with second-order operators, third-order convergence with third-order operators, fourth-order convergence with  fourth-order  operators and fifth-order convergence with fifth-order operators.

The CPU times associated with each computation are reported in Table~\ref{tab:CPUNonCompactJunction}. The CPU times corresponding to a relative error of the order of $10^{-6}$ are also highlighted in bold font. For the same accuracy, the CPU time associated with the second-order operators is 2.5 times larger than the ones associated with the higher-order operators.

%
%
%

\begin{table}[htdp]
\begin{center}
\begin{tabular}{|c|c|c|c|c|}
\hline
 N& order 2 & order 3 & order 4 & order 5 \\
 \hline
 32 &0.707003  &  1.359596  &2.381819&  3.028414 \\
64 &0.969357  &   2.641288 &  3.757290 &  5.092958\\
128  & 2.135884&  3.319466 & 3.876460 &  4.843739\\
 256 & 2.099840 &   3.097754 & 4.138360&4.951223\\
 512 & 2.051456  &  3.043053  & 4.327195   & 5.108929 \\
 \hline
 \end{tabular} 
 \caption{Numerical convergence behavior for each operator choice}
\label{tab:CVorderNonCompactJunction}
\end{center}
\end{table}

\begin{table}[htdp]
\begin{center}
\begin{tabular}{|c|c|c|c|c|}
\hline
 N& order 2 & order 3 & order 4 & order 5 \\
 \hline
 16 & 33.05 & 34.30 & 34.26 & 34.99 \\
32 & 33.01& 34.54 &  35.71& {\bf 35.74} \\
64  & 34.88 &  35.84&  {\bf 38.97}&78.98 \\
 128 &   35.78 & {\bf 39.48} & 40.22  &  86.54\\
 256 &  39.47 &  45.09 &  47.87 &111.73\\
  512 & {\bf 93.04} &  125.99 &   126.18 & 190.01\\
 \hline
 \end{tabular} 
 \caption{CPU times (in seconds) for each operator choice in the junction problem}\label{tab:CPUNonCompactJunction}
\end{center}
\end{table}
\subsection{Cable and soma}
\begin{figure}
\begin{center}
\includegraphics[height=3.0in]{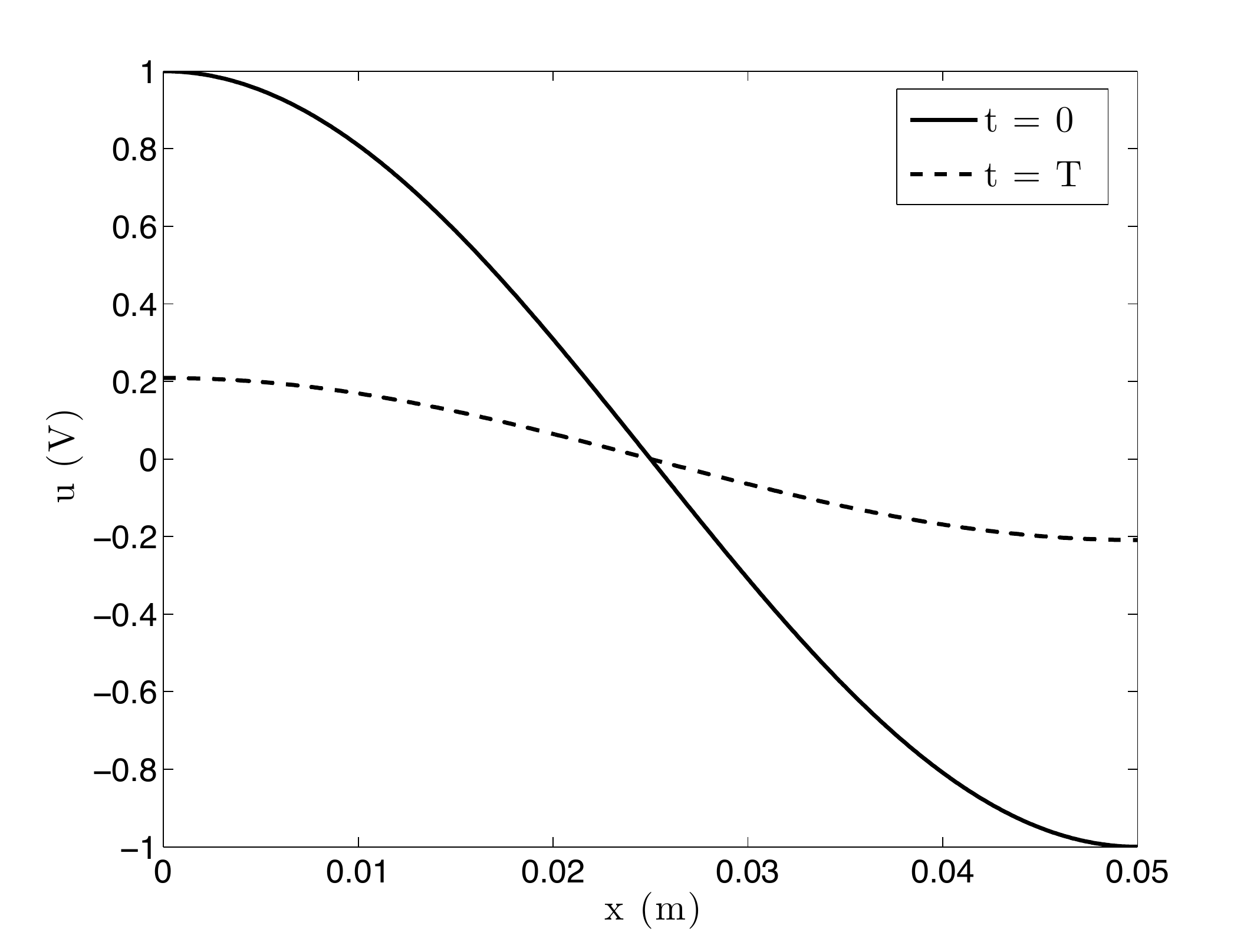}
 \caption{Exact solutions for the cable with soma problem}
 \label{fig:exactsolCableSoma}
 \end{center}
\end{figure}

The exact manufactured solution $u(x,t)$ and the corresponding initial boundary value problem are specified in Appendix D.  The physical constants associated with the problem are reported in Tables~\ref{tab:physvar} and~\ref{tab:cablesomaGeom}. The exact solutions $u(x,t)$ at times $t=0$ and $t=T$ are depicted in Figure~\ref{fig:exactsolCableSoma}.

\begin{table}[htdp]
\begin{center}\begin{tabular}{|c|c|c|}
\hline
Cable length & $L$ & $0.05$ m \\
\hline
Cable radius & $a(x)$ & $0.476\times 10^{-3}$ m  \\
\hline
Soma radius & $r_{\text{soma}}$ & $2\times 10^{-3}$ m \\
\hline
\end{tabular} 
\end{center}\caption{Geometric parameters associated with the cable with soma problem}
\label{tab:cablesomaGeom}
\end{table}

\begin{figure}
\begin{center}
\includegraphics[height=3.0in]{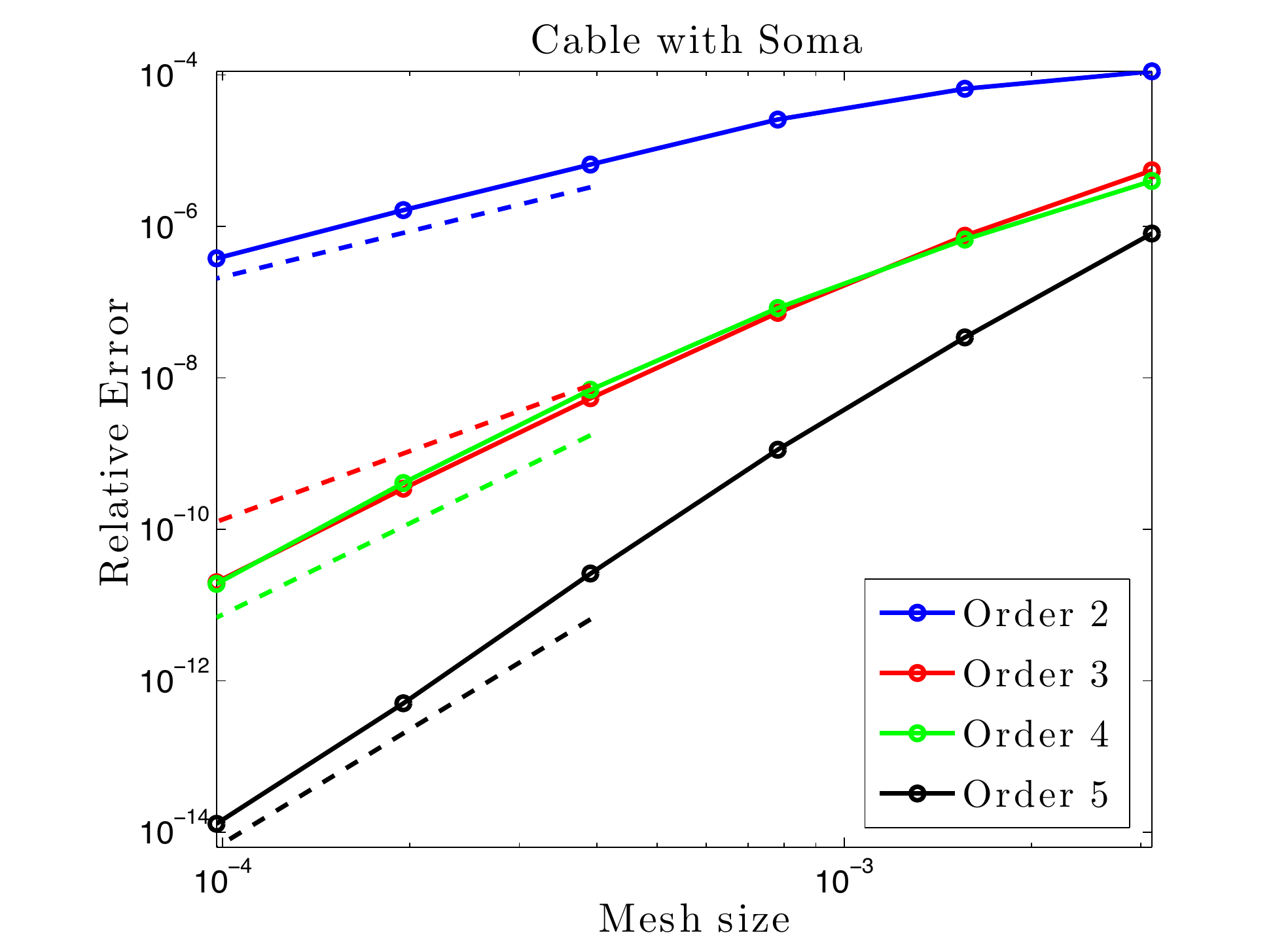}
 \caption{Errors in function of mesh size and operators order for the cable attached to soma problem}
 \label{fig:CableSomaCVOrderNonCompact}
 \end{center}
\end{figure}

The convergence results are reported in Table~\ref{tab:CVorderNonCompactCableSoma} and Figure~\ref{fig:CableSomaCVOrderNonCompact}, in which slopes corresponding to the operators orders are also depicted.  Second-order convergence can be observed using second-order operators, fourth order convergence is observed when using third- and fourth-order operators and fifth-order convergence for fifth-order operators. For this case, no real advantage can be seen in terms of using fourth-order operators over third-order operators.

The CPU times associated with each computation are reported in Table~\ref{tab:CPUNonCompactCableSoma}. In that table, the CPU times corresponding to a relative error of the order of $5\times 10^{-7}$ are highlighted in bold font. Note that for the same accuracy, the CPU time associated with the second-order operators is 20 times larger than the ones associated with the high-order operators.

%
%

\begin{table}[htdp]
\begin{center}
\begin{tabular}{|c|c|c|c|c|}
\hline
 N& order 2 & order 3 & order 4 & order 5 \\
 \hline
32 & 0.766683& 2.862748 & 2.567863 & 4.554248 \\
64  & 1.345940 &  3.386744&  3.004467 &4.927776 \\
 128 &  1.978948& 3.752389 & 3.588062  &  5.442091\\
 256 &  2.005454 &  3.962842 &  4.090079 &5.691623\\
  512 &2.000201 &  3.983892 &   4.314146 & 5.161448\\
 \hline
 \end{tabular} 
 \caption{Numerical convergence behavior for each operator choice}\label{tab:CVorderNonCompactCableSoma}
\end{center}
\end{table}

\begin{table}[htdp]
\begin{center}
\begin{tabular}{|c|c|c|c|c|}
\hline
 N& order 2 & order 3 & order 4 & order 5 \\
 \hline
 16 & 17.03 & 16.64 & 16.03 & {\bf 16.27} \\
32 & 16.81& {\bf 16.78} & {\bf 16.50} & 16.67 \\
64  & 20.25 &  21.49&  20.47&36.76 \\
 128 &  24.05 & 24.80 & 29.57  &  42.30\\
 256 &  39.76 &  44.88 &  49.14 &117.77\\
  512 & {\bf 315.09} &  294.97 &   318.83 & 551.21\\
 \hline
 \end{tabular} 
 \caption{CPU times (in seconds) for each operator choice in the cable attached to soma problem}\label{tab:CPUNonCompactCableSoma}
\end{center}
\end{table}

\section{Numerical applications}\label{sec:Appli}

\subsection{Time discretization}
After  space discretization by the approach developed in this paper, the Hodgkin-Huxley equations become a set of coupled ODEs of the form
\begin{gather*}
\begin{split}
\frac{d \ubold}{dt} &= \Abold_1(\wbold)\ubold + \bbold_1(\wbold,t) \\
\frac{d \wbold}{dt} &= \Abold_2(\ubold)\wbold + \bbold_2(\ubold) 
\end{split}
\end{gather*}
where the discrete variables for the potential and gating variables are $\ubold\in\mathbb{R}^{N+1}$ and $\wbold = [\mbold,\hbold,\nbold]^T\in\mathbb{R}^{3(N+1)}$, respectively.

This set of coupled ODEs is then discretized in time, using the second-order staggered difference scheme initially developed by Hines~\cite{hines84}. The scheme is
\begin{gather*}
\begin{split}
\left(\Ibold - \frac{\Delta t}{2} \Abold_2(\ubold^{n})\right) \wbold^{n+\frac{1}{2}} &= \left(\Ibold +\frac{\Delta t}{2}\Abold_2(\ubold^{n})\right) \wbold^{n-\frac{1}{2}} + \Delta t \bbold_2(\ubold^{n})   \\
\left(\Ibold - \frac{\Delta t}{2} \Abold_1(\wbold^{n+\frac{1}{2}})\right) \ubold^{n+1} &= \left(\Ibold +\frac{\Delta t}{2}\Abold_1(\wbold^{n+\frac{1}{2}})\right) \ubold^{n} + \Delta t \bbold_1(\wbold^{n+\frac{1}{2}},t^{n+\frac{1}{2}}).
\end{split}
\end{gather*}

Here $\displaystyle{\ubold^n \approx \ubold(t^n) = \ubold(n\Delta t)}$ and $\displaystyle{\wbold^{n+\frac{1}{2}} \approx \wbold(t^{n+\frac{1}{2}}) = \wbold\left(\left(n+\frac{1}{2}\right)\Delta t\right)}$.


\subsection{Potential propagation in a cable attached to the soma}

The physical constants associated with this problem are provided in Table~\ref{tab:physvar}. These correspond to the giant squid's axon studied by Hodgkin and Huxley in~\cite{hodgkin52}. The expressions for the gating functions are provided in Appendix A.

\begin{figure}
\begin{center}
\includegraphics[height=3.5in]{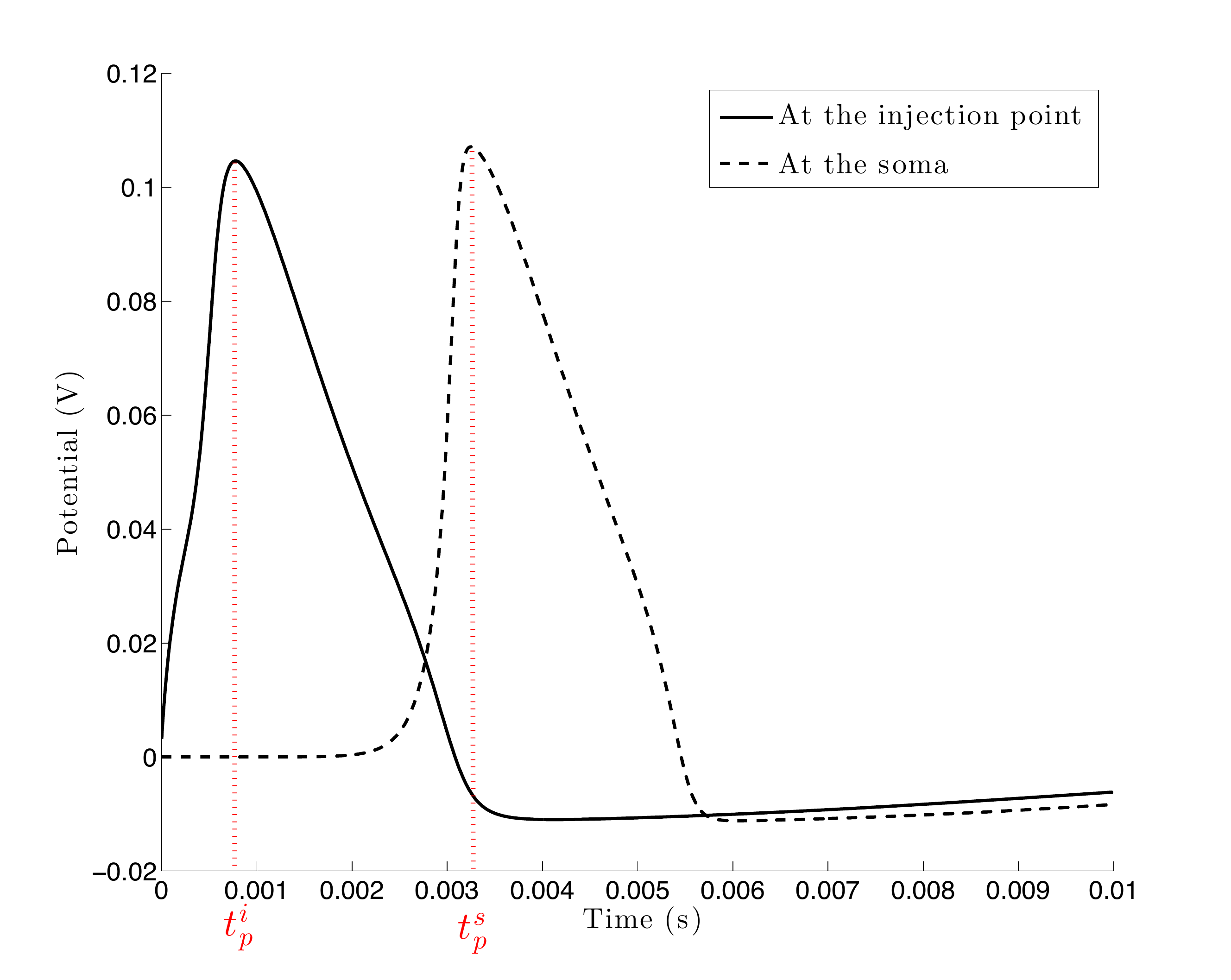}
 \caption{Action potential at the current injection point and the soma obtained using fifth-order operators}
 \label{fig:cableSomaTH}
 \end{center}
\end{figure}

Current with Gaussian distribution and maximum intensity $I = 10$ nA is injected at the top extremity of the cable depicted in Figure~\ref{fig:Configurations}(a) and the cable potential time history for $t = [0,10]$ ms is recorded at the two extremities of the cable. The corresponding time histories obtained by applying the procedure developed in this paper with $N=32$ and 5-th order differential operators and a time step size of $dt = 10^{-2}$ ms are reported in Figure~\ref{fig:cableSomaTH}. The typical action potential characteristic shape~\cite{keener09} can be observed.

\begin{figure}
\begin{center}
\includegraphics[height=3.0in]{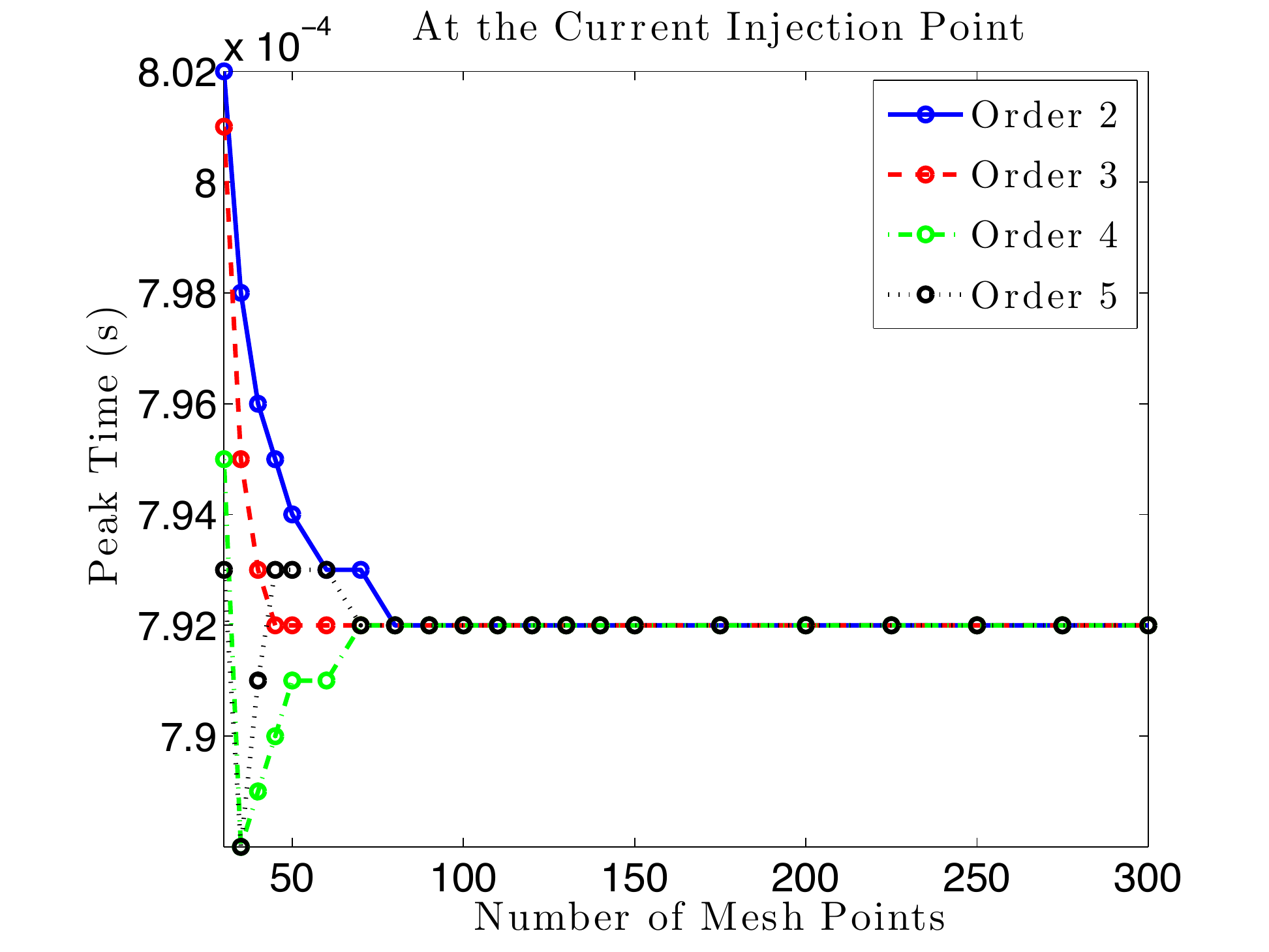}
 \caption{Action potential peak times $t_p^i$ at the current injection point}
 \label{fig:injectionPeakTimes}
 \end{center}
\end{figure}
\begin{figure}
\begin{center}
\includegraphics[height=3.0in]{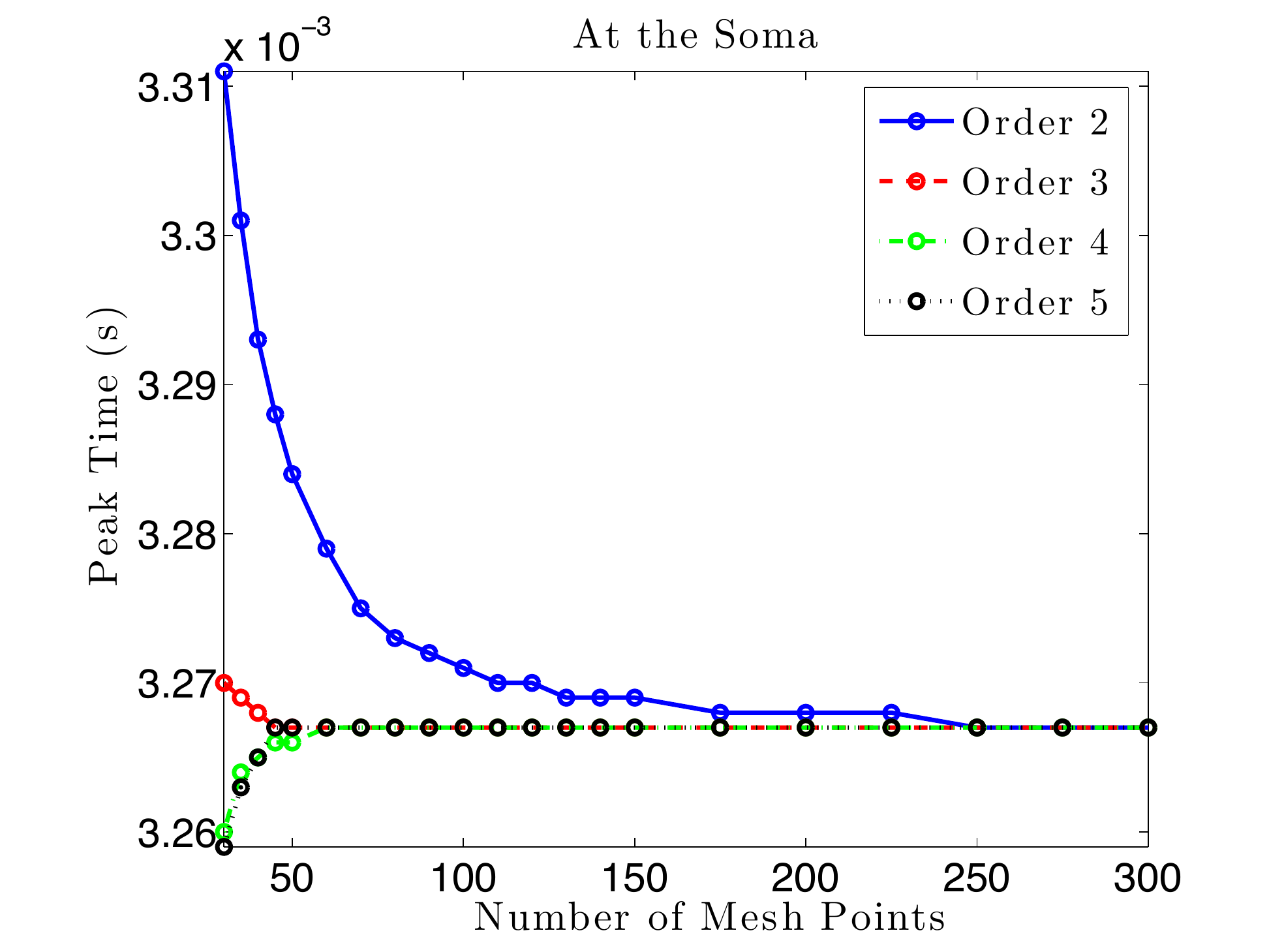}
 \caption{Action potential peak times $t_p^s$ at the soma}
 \label{fig:somaPeakTimes}
 \end{center}
\end{figure}
\begin{figure}
\begin{center}
\includegraphics[height=3.0in]{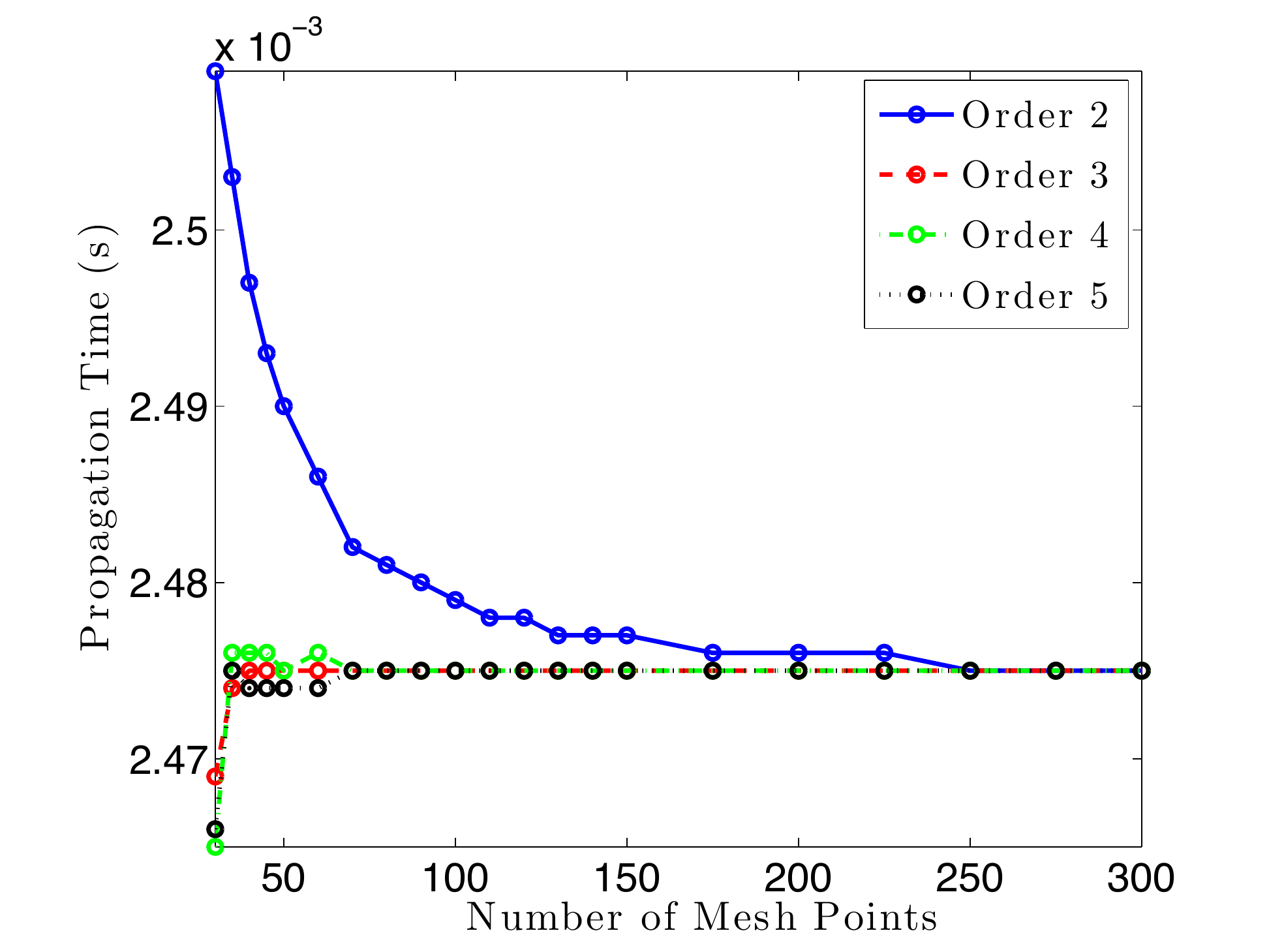}
 \caption{Propagation times from the current injection point to  the soma}
 \label{fig:propagationTimes}
 \end{center}
\end{figure}

Next, the action potential peak times obtained at the two extremities of the cable of interest  are recorded and shown  for  various choices of order for the differential operators and various numbers of mesh points. The peak time at the injection points $t_p^i$,  depicted in Figure~\ref{fig:cableSomaTH}, are almost identical as it can be observed in Figure~\ref{fig:injectionPeakTimes}. However, for coarse meshes, second-order operators result in an inaccurate peak time $t^s_p$ at the Soma, as observed in Figure~\ref{fig:somaPeakTimes}. This results in inaccurate propagation times for these coarse meshes when second-order operators are used, whereas the propagations times are accurate with high-order operators even when using those coarse meshes, as shown in Figure~\ref{fig:propagationTimes}.

\subsection{Potential propagation in a dendritic network}

\begin{figure}
\begin{center}
\includegraphics[height=3.0in]{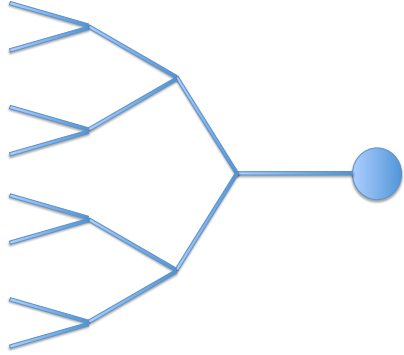}
 \caption{Dendritic tree connected to the soma}
 \label{fig:large_tree}
 \end{center}
\end{figure}

In order to illustrate the capability of the proposed approach to predict potential propagation in large dendritic networks, a tree with 15 branches attached to the soma, as depicted in Figure~\ref{fig:large_tree}, is considered. The dimensions of each branch in the tree are chosen according to Rall's 3/2 power law~\cite{rall60} and correspond to Rallpack's second configuration~\cite{bhalla92}. The cable lengths and radii at each level of the tree are reported in Table~\ref{tab:treebranchesdim}. All the other physical constants are identical as in the previous numerical experiments. Synaptic current input is modeled by a current injection of amplitude $I=500$ nA for a duration of $1$ ms at the extremity of each of the 8 extremal branches. Each branch is discretized using 31 nodes, resulting in 1860 degrees of freedom. Fifth order difference operators are here used. The potential propagation is subsequently simulated in the entire tree in the time interval for $t\in[0,0.2]$ s using a time step $\Delta t = 10^{-4}$ s. The entire simulation takes about 14 s CPU time when performed with MATLAB running on a Mac Book Pro 2.9 GHz Intel Core i7, 8GB 1600 MHz DDR3. 

\begin{figure}
\begin{center}
 \begin{subfigmatrix}{2}
  \subfigure[Current input activity for each of the extremal branches (each line segment is of length $1$ ms)]{\includegraphics[height=2.5in]{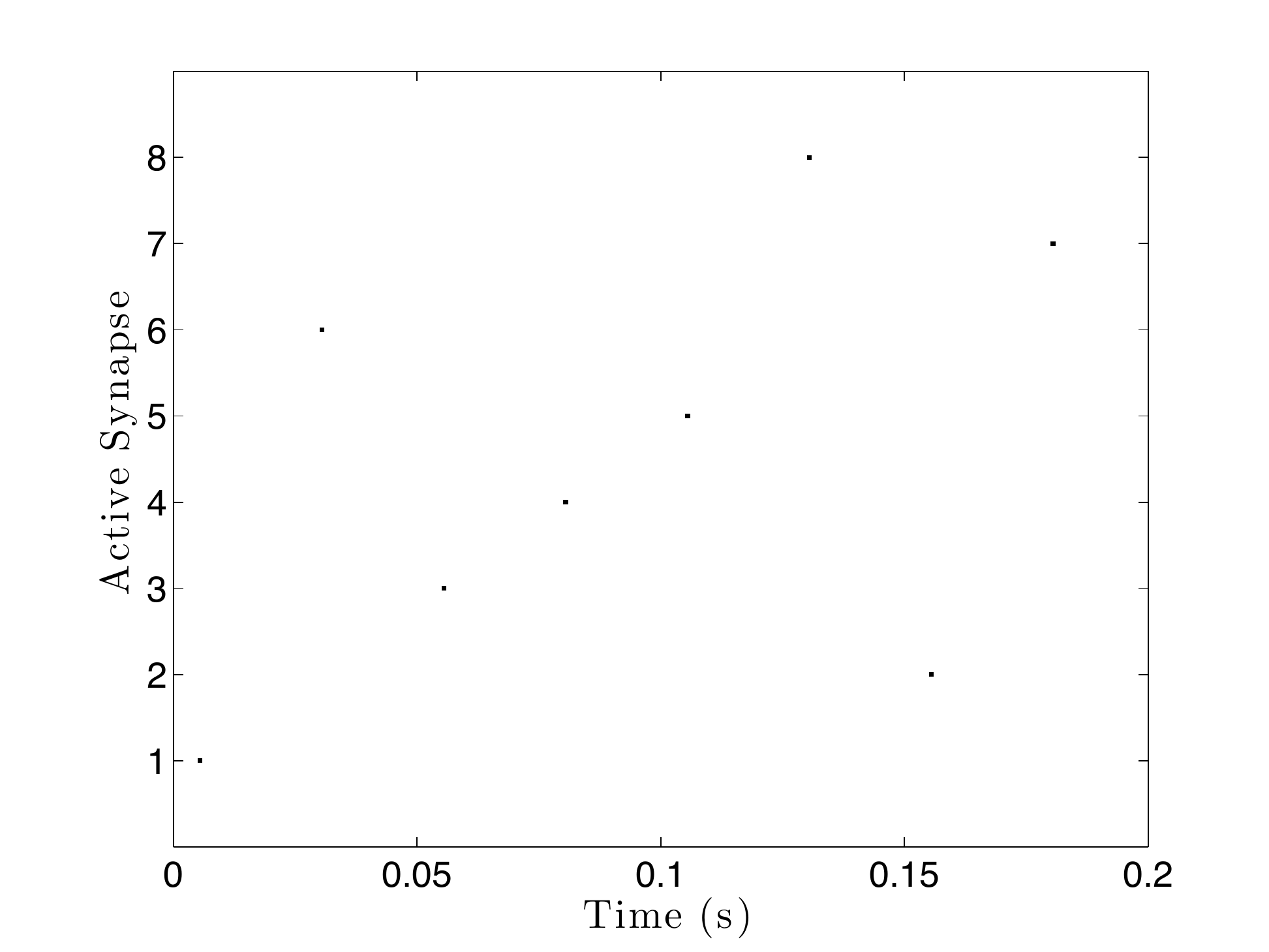}}
  \subfigure[Action potential at the soma]{ \includegraphics[height=2.5in]{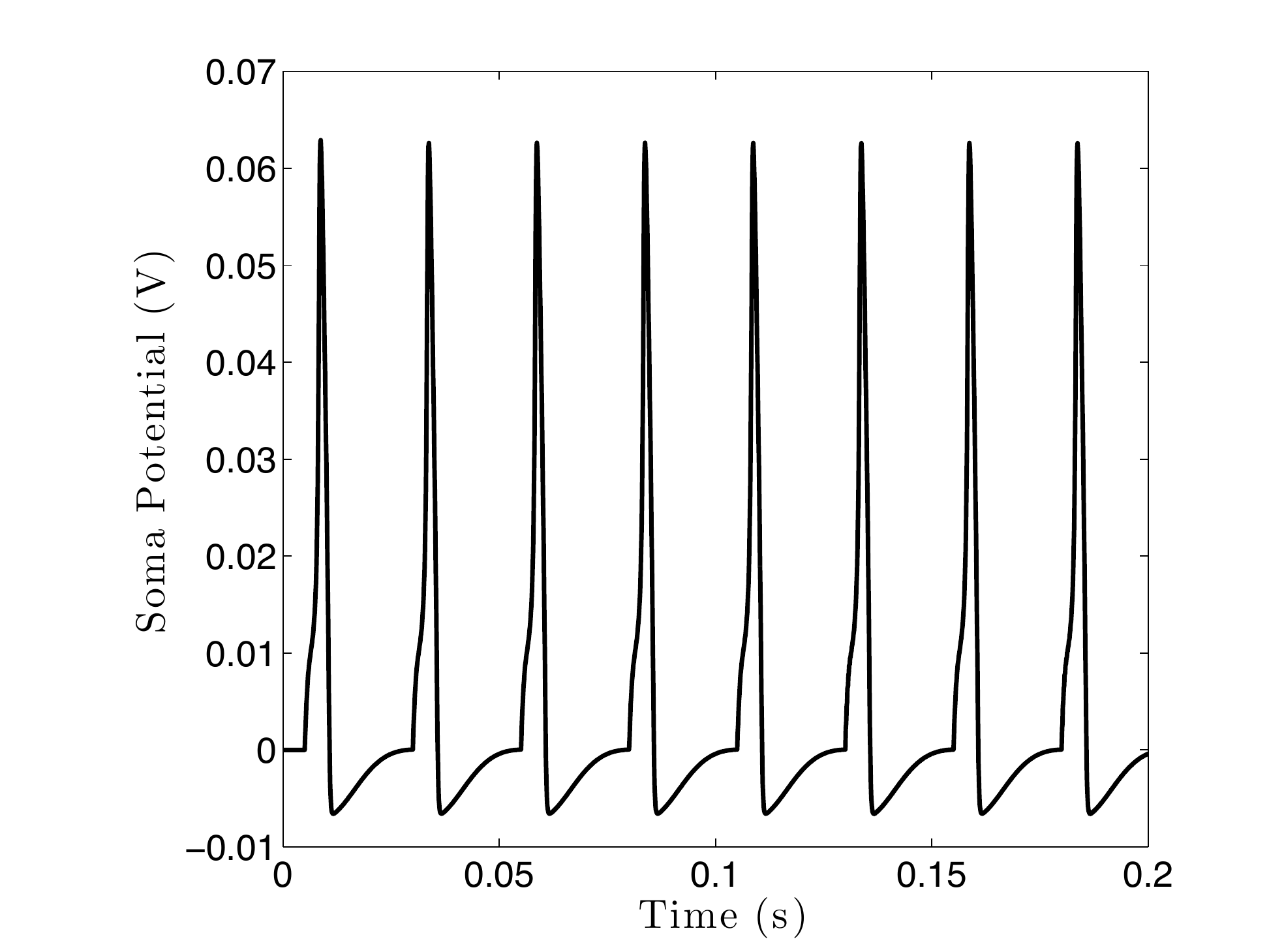}}
  \end{subfigmatrix}
 \caption{Potential propagation in the network for the first current input configuration}
 \label{fig:Configuration1}
 \end{center}
\end{figure}

Two current input configurations are considered. The first one is depicted in Figure~\ref{fig:Configuration1} (a). In that figure, the interval of duration $1$ ms, at which a given branch input is active, is represented by a line segment. The corresponding potential recorded at the soma  is reported in Figure~\ref{fig:Configuration1} (b). There are eight characteristic action potential spikes, each one corresponding to a specific extremal branch input.

\begin{figure}
\begin{center}
 \begin{subfigmatrix}{2}
  \subfigure[Current input activity for each of the extremal branches (each line segment is of length $1$ ms)]{\includegraphics[height=2.5in]{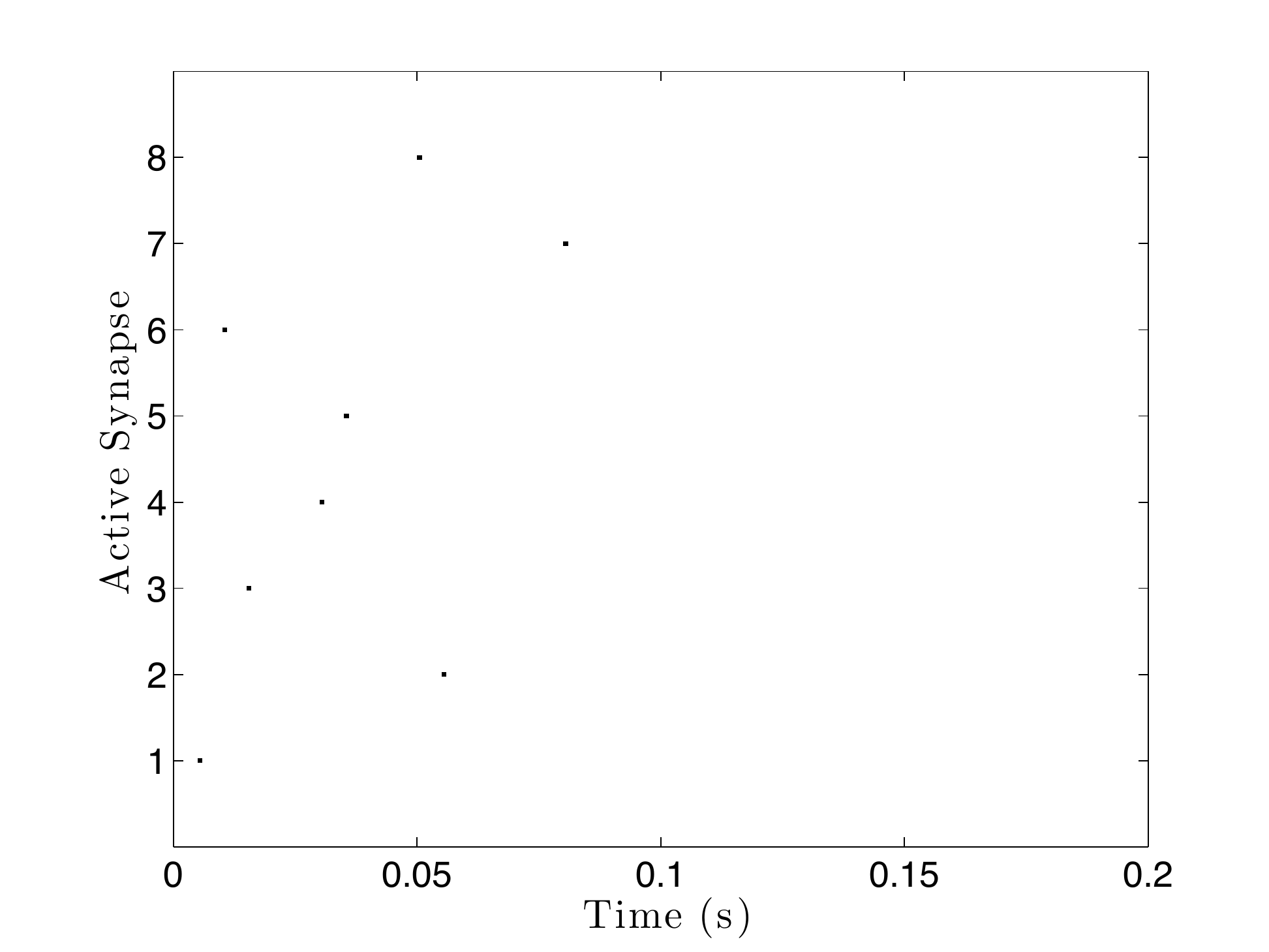}}
  \subfigure[Action potential at the soma]{ \includegraphics[height=2.5in]{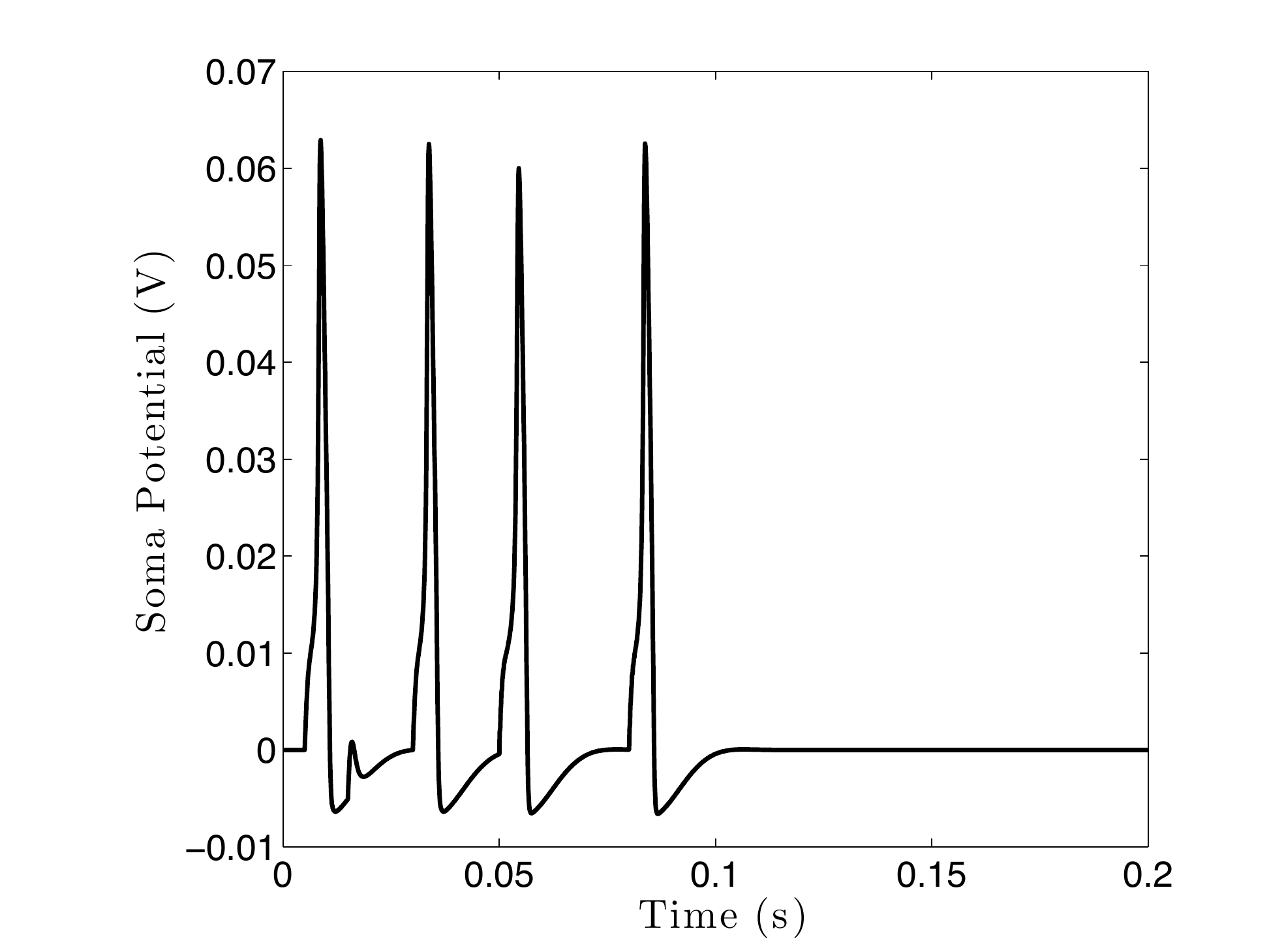}}
  \end{subfigmatrix}
 \caption{Potential propagation in the network for the second current input configuration}
 \label{fig:Configuration2}
 \end{center}
\end{figure}

The second numerical experiment is carried out with the  current input activity reported in Figure~\ref{fig:Configuration2} (a).  In this case, the current input is limited to the time interval $t\in[0,0.1]$ s. The associated potential at the soma is presented in Figure~\ref{fig:Configuration2} (b). Only four characteristic action potential spikes are propagated from the extremity branches of the tree to the soma. This numerical experimental emphasizes the nonlinear nature of the Hodgkin-Huxley equations. It also illustrates the concept of refractory period for action potentials observed in neurons~\cite{koch04}: after an action potential has been generated, the membrane is hyper-polarized and there is a minimal period of time after which another action potential can be created. In the second numerical experiment, although eight current inputs were injected, four failed to create full-amplitude spikes since their respective injection times were within the refractory periods of the four action potentials. A small spike can also be observed around $t=0.017$ s, but it is not a full action potential. In the first experiment, however, the time between two current injections is longer, leading to eight spikes.

\begin{table}[htdp]
\begin{center}
\begin{tabular}{|c|c|c|c|}
\hline
 Level& Number of branches & Branch length (m) & Branch radius (m)\\
 \hline
 1 & 1 & $32.0\times 10^{-6}$ &$8.0\times 10^{-6}$  \\
2 & 2 &$25.4\times 10^{-6}$& $5.04\times 10^{-6}$ \\
3  & 4 & $20.16\times 10^{-6}$ &  $3.18\times 10^{-6}$  \\
 4 &  8 & $16.0 \times 10^{-6}$ & $2.0\times 10^{-6}$\\
 \hline
 \end{tabular} 
 \caption{Dimensions of each branch in the dendritic tree}\label{tab:treebranchesdim}
\end{center}
\end{table}

\section{Conclusions}\label{sec:Concl}

Stable and accurate high-order schemes based on Summation-By-Parts (SBP) form with the weak Simultaneous Approximation Term (SAT) procedure for boundary and interface conditions are developed for the Hodgkin-Huxley equations.  These schemes are shown in this paper to be easily applicable to neuronal systems constituted of multiple branches such as in the case of dendritic trees. Well-posedness of the continuous problem as well as energy estimates are established for dendritic trees connected  to the soma.  The SBP-SAT procedure leads to an energy estimate and proof of stability of the discrete scheme. 

Convergence and order of accuracy are illustrated on model problems. Important gains in CPU times for a given level of accuracy are also demonstrated when high-order operators are used. The superiority of high-order operators is also illustrated on a simple problem constituted of a cable connected to the soma. Finally, the proposed approach is applied to a realistic dendritic tree configuration constituted of 15 cables connected to the soma. On that configuration, the nonlinearity properties of the Hodgkin-Huxley equations are numerically demonstrated.

%
\section*{Appendix A: expressions for the gating functions}
%
  The expressions for the coefficients $\alpha_m$, $\alpha_h$, $\alpha_n$, $\beta_m$, $\beta_h$ and $\beta_n$ were determined for the giant squid axon by Hodgkin and Huxley~\cite{hodgkin52} as:
\begin{gather}
\begin{split}
\alpha_m(u)  &= \frac{10^5(0.025-u)}{e^{\frac{0.025-u}{0.01}}-1},~~\beta_m(u) = 4\times10^3e^{-\frac{u}{0.018}} \\
\alpha_h(u)  &= 70e^{-\frac{u}{0.02}},~~\beta_h(u) \frac{10^3}{e^{\frac{0.03-u}{0.01}}+1} \\
\alpha_n(u)  &= \frac{10^4(0.01-u)}{e^{\frac{0.01-u}{0.010}}-1},~~\beta_n(u) = 125 e^{-\frac{u}{0.08}}. 
  \end{split}
  \end{gather}
\section*{Appendix B: proof of proposition 1}

The differential equations associated with $m(x,t)$ are a set of uncoupled ODEs for each $x$. As a result, it is sufficient to consider one of those ODE only for $x_0\in[0,L]$:
\begin{equation}
 m_t(x_0,t) =\alpha_m(u(x_0,t)) (1-m(x_0,t)) - \beta_m(u(x_0,t))m(x_0,t) , ~~t\in[0,T].
 \end{equation}

 To simplify the notations, in the rest of this section, the dependency on $x_0$ is dropped: $\alpha_m(u(x_0,t)) = \alpha_m(u(t))$,  $\beta_m(u(x_0,t)) = \beta_m(u(t))$ and $ m(x_0,t) = m(t)$.

From the expressions of $\alpha_m$ and $\beta_m$ given in Appendix 1, it can be shown  that
\begin{equation}
 \alpha_m(0) > 0,~~\beta_m(1) > 0.
 \end{equation}
These properties also hold  for $\alpha_h$, $\alpha_n$, $\beta_h$ and $\beta_n$.

An initial condition $ m(0) \in [0,1]$ is then considered.  Assume that there exists $t_1\in(0,T]$ such that $m(t_1)>1$. As $m(\cdot)$ is continuous, this means that there exists $0<t_0<t_1<T$ such that $m(t_0)=1$ and $m(t_0+\epsilon)>1,~\forall 0<\epsilon\leq t_1-t_0$. Since
\begin{equation}
 \frac{d m}{dt}(t_0) = \alpha_m(u(t_0))(1-1)-\beta_m(u(t_0))1 =  -\beta_m(u(t_0))< 0,
 \end{equation}
this is a contradiction. Hence $m(t)\leq 1,~\forall t\in[0,T]$.

 Similarly, assume now  that there exists $t_1$ such that $m(t_1)<0$. As $m(t)$ is continuous, that means that there exists $0<t_0<t_1$ such that $m(t_0)=0$ and $m(t_0+\epsilon)<0,~\forall 0<\epsilon\leq t_1-t_0$. However,
\[ \frac{d m}{dt}(t_0) = \alpha_m(u(t_0))(1-0) - \beta_m(u(t_0))0 = \alpha_m(u(t_0)) > 0\]
which leads to  a contradiction as well.

 In conclusion $0 \leq m(t) \leq 1,~\forall t\in[0,T]$. This property also holds for $h(t)$ and $n(t)$, respectively.

\section*{Appendix C: manufactured solution for the junction problem}
A case with 3 branches having a junction at $x^{(i)}=0,~i=1,\cdots,3$ is considered. The geometric parameters of the problem are  $a^{(1)}(x^{(1)})= a^{(2)} (x^{(2)})= 2^{-2/3} a^{(3)}(x^{(3)}) = a_0$, where $x^{(i)}\in[0,L^{(i)}]$, and their lengths are given by $L^{(i)} = L,~i=1,2$ and $L^{(3)} = 2^{1/3} L$. The following manufactured solutions are considered, for a problem  without current injection, that is $I(t)=0,~t\in[0,T]$, where $T=10^{-5}$ s:
\begin{gather}\label{eq:exactsolJunction}
\begin{split}
u^{(i)}(x^{(i)},t) &= \text{exp}\left[-t\left( \frac{1}{C_m}\sum_{i=j}^3g_j+\mu a_0\left(\frac{3\pi}{2L}\right)^2\right)\right] \sin\left(\frac{3\pi x^{(i)}}{2L}\right),~~i=1,2 \\
u^{(3)}(x^{(3)},t) &= \text{exp}\left[-t\left( \frac{1}{C_m}\sum_{i=j}^3g_j+ \mu a_0\left(\frac{3\pi}{2L}\right)^2\right)\right] \sin\left(-\frac{3\pi x^{(3)}}{2^{4/3}L}\right),\\
m^{(i)}(x^{(i)},t) &= 1,~~i=1,\cdots,3,\\
h^{(i)}(x^{(i)},t) &= 1,~~i=1,\cdots,3,\\
n^{(i)}(x^{(i)},t) &= 1,~~i=1,\cdots,3,
\end{split}
\end{gather}
with $(x^{(i)},t) \in [0,L^{(i)}]\times [0,T],~i=1,\cdots,3$. 

The solutions~(\ref{eq:exactsolJunction}) satisfy the equations
\begin{gather}
\begin{split}
u^{(i)}_t &=  \frac{\mu}{a^{(i)}(x^{(i)})} \left(a^{(i)}(x^{(i)})^2u^{(i)}_{x^{(i)}}\right)_{x^{(i)}} - \frac{1}{C_m}g(m^{(i)},h^{(i)},n^{(i)}) u^{(i)} + f(m^{(i)},h^{(i)},n^{(i)},t) + F_u \\
m^{(i)}_t &= \alpha_m(u^{(i)})(1-m^{(i)}) - \beta_m(u^{(i)}) m^{(i)} + F^{(i)}_m(u^{(i)}) \\
h^{(i)}_t &= \alpha_h(u^{(i)})(1-h^{(i)}) - \beta_h(u^{(i)}) h^{(i)} + F^{(i)}_h(u^{(i)}) \\
n^{(i)}_t &= \alpha_n(u^{(i)})(1-n^{(i)}) - \beta_n(u^{(i)}) n^{(i)} + F^{(i)}_n(u^{(i)}), \\
\end{split}
\end{gather}
where 
\begin{gather}
\begin{split}
F_u &= - \frac{1}{C_m} \sum_{j=1}^3 g_j E_j \\
F^{(i)}_m(t) &= \beta^{(i)}_m(u^{(i)}) \\
F^{(i)}_h(t) &= \beta^{(i)}_h(u^{(i)}) \\
F^{(i)}_n(t) &= \beta^{(i)}_n(u^{(i)}) \\
\end{split}
\end{gather}
and the boundary conditions are
\begin{gather}
\begin{split}
u^{(i)}(0,t) &= 0,\\
\sum_{i=1}^3 \left(a^{(i)}\right)^2 u^{(i)}_{x^{(i)}}(0) &=0, \\
u^{(i)}_{x^{(i)}}(L^{(i)}) &= 0,~i=1,2, \\
u^{(3)} (x^{(3)},t) &=  exp\left[-t\left( \frac{1}{C_m}\sum_{i=j}^3g_j+ \mu a_0\left(\frac{3\pi}{2L}\right)^2\right)\right]  .
\end{split}
\end{gather}
The numerical values chosen for each physical constant are the same as the ones reported in Table~\ref{tab:physvar} and~\ref{tab:junctionGeom}.
\section*{Appendix D: manufactured solution for the cable with soma problem}

A problem with constant cable radius $a(x) = a_0,~x\in[0,L]$ and without current injection, that is $I(t)=0,~t\in[0,T]$, where $T=10^{-5}$ s:
\begin{gather}\label{eq:exactsolCableSoma}
\begin{split}
u(x,t) &= \text{exp}\left[-t\left( \frac{1}{C_m}\sum_{j=1}^3g_j+\mu a_0 \left(\frac{\beta}{L}\right)^2\right)\right] \cos\left(\frac{\beta x}{L}\right) \\
m(x,t) &= 1\\
h(x,t) &= 1\\
n(x,t) &= 1
\end{split}
\end{gather}
where $(x,t)\in [0,L]\times [0,T]$ and $\beta$ satisfies $\frac{\tan \beta}{\beta} = - \frac{\mu}{\eta a_0 L}$.

The solutions~(\ref{eq:exactsolCableSoma}) satisfy the equations
\begin{gather}
\begin{split}
u_t &=  \frac{\mu}{a(x)} \left(a(x)^2u_x\right)_{x} -\frac{1}{C_m} g(m,h,n) u + f(m,h,n,x, t) + F_u(t) \\
m_t &= \alpha_m(u)(1-m) - \beta_m(u) m + F_m(u) \\
h_t &= \alpha_h(u)(1-h) - \beta_h(u) h + F_h(u) \\
n_t &= \alpha_n(u)(1-n) - \beta_n(u) n + F_n(u) \\
\end{split}
\end{gather}
where 
\begin{gather}
\begin{split}
F_u(t) &= -  \frac{1}{C_m}\sum_{j=1}^3 g_j E_j \\
F_m(t) &= \beta_m(u) \\
F_h(t) &= \beta_h(u) \\
F_n(t) &= \beta_n(u) \\
\end{split}
\end{gather}
and the boundary conditions
\begin{gather}
\begin{split}
u_x(0) &= 0\\
u_t(L) &= -\eta a(L)^2 u_x(L) -\frac{1}{C_m} g(m(L),h(L),n(L)) + f(m(L),h(L),n(L),L,t) + F_b(t)
\end{split}
\end{gather}
where 
\begin{equation}
F_b(t) = -\frac{1}{C_m}\sum_{j=1}^3 g_j E_j.
\end{equation}
The numerical values chosen for each physical constant are provided in Tables~\ref{tab:physvar} and~\ref{tab:cablesomaGeom}.

\bibliography{pap}
\bibliographystyle{model1-num-names.bst}

\end{document}